\documentclass[11pt]{amsart}
\usepackage{amscd,amssymb,verbatim}

\theoremstyle{plain}
\newtheorem{Thm}{Theorem}[section]
\newtheorem{Cor}[Thm]{Corollary}
\newtheorem{Lem}[Thm]{Lemma}
\newtheorem{Prop}[Thm]{Proposition}
\newtheorem{Def}[Thm]{Definition}

\newtheorem{remark}{Remark}

\begin{document}

\title{On the complemented subspaces of the Schreier spaces}
\author{I. Gasparis}
\address{Department of Mathematics \\
         Oklahoma State University \\
         Stillwater, OK 74078-1058 \\
         U.S.A.}
\email{ioagaspa@math.okstate.edu}
\author{D. H. Leung}
\address{Department of Mathematics \\
         National University of Singapore \\
         Singapore 117543}
\email{matlhh@nus.edu.sg}
\date{July 3, 1999}
\begin{abstract}
It is shown that for every \(1 \leq \xi < \omega\) the Schreier space 
\(X^\xi\) admits a set of continuum cardinality whose elements are
mutually incomparable
complemented subspaces spanned by subsequences of \((e_{n}^{\xi})\), the 
natural Schauder basis of \(X^\xi\). It is also shown that there exists
a complemented subspace spanned by a block basis of \((e_{n}^{\xi})\),
which is not isomorphic to a subspace generated by a subsequence 
of \((e_n^\zeta)\), for every
\(0 \leq \zeta \leq \xi\).
Finally, an example is given of an uncomplemented
subspace of \(X^\xi\) which is spanned by a block basis of \((e_{n}^{\xi})\).
\end{abstract}
\maketitle
\section{Introduction} \label{S:1}
The Schreier families \(\{S_\xi\}_{\xi < \omega_1}\) of finite
subsets of positive integers (the precise definition is given in
the next section), introduced in \cite{aa}, have 
played a central role in the development of modern Banach space
theory. We mention the use of Schreier families in the construction
of mixed Tsirelson spaces which are asymptotic \(\ell_1\) and arbitrarily
distortable \cite{ad}. The distortion of mixed Tsirelson spaces has
been extensively studied in \cite{ano}. In that paper as well as in
\cite{otw}, the moduli \((\delta_\alpha)_{\alpha < \omega_1}\) were
introduced measuring the complexity of the asymptotic \(\ell_1\) structure
of a Banach space. The definitions of those moduli also involve the
Schreier families. Other applications can be found in \cite{amt} and
\cite{ag} where the Schreier families form the main tool for 
determining the structure of those convex combinations of a weakly
null sequence that tend to zero in norm, or are equivalent to the
unit vector basis of \(c_0\).
For applications of the Schreier families in the construction of
hereditarily indecomposable Banach spaces, we refer to \cite{ad}
and \cite{af}. 

A notion companion to the Schreier families is that of the 
Schreier spaces. These are Banach spaces whose norm is related
to a corresponding Schreier family. More precisely, for every
countable ordinal \(\xi\), we define a norm \(\|\cdot\|_\xi\)
on \(c_{00}\), the space of finitely supported real valued sequences,
in the following manner: Given \(x=\bigl(x(n)\bigr) \in c_{00}\) define
\[ \|x\|_\xi = \sup_{F \in S_\xi} \sum_{n \in F} |x(n)|.\]
\(X^\xi\), the Screier space of order \(\xi\), is the completion
of \(c_{00}\) under the norm \(\|\cdot\|_\xi\).
\(X^0 = c_0\), the Banach space of null sequences.
\(X^1\) was first considered by Schreier \cite{sc} in order
to provide an example of a weakly null sequence without
Cesaro summable subsequence. It is proven in \cite{aa} that the
natural Schauder basis \((e_{n}^{\xi})\) of \(X^\xi\) is  
1-unconditional and shrinking.
\(X^1\) has been studied in \cite{o2} where it is shown that
every quotient of \(X^1\) is \(c_0\)-saturated. That is, every
infinite dimensional subspace contains a further subspace isomorphic
to \(c_0\).

Given \(M\), an infinite subset of
 \(\mathbb{N}\), we let \(X_{M}^{\xi}\) denote the
closed linear subspace of \(X^\xi\) spanned by the subsequence
\((e_{n}^{\xi})_{n \in M}\). For an element \(x \in X^{\xi}\),
\(x= \sum_{n \in \mathbb{N}} a_n e_{n}^{\xi}\), we set
\(\|x\|_0 = \sup_{n \in \mathbb{N}} |a_n|\).
The main result of this paper is
the following 
\begin{Thm} \label{T:1}
Let \(L=(l_n)\), \(M=(m_n)\) be infinite subsets of \(\mathbb{N}\),
and let \(\xi < \omega\).
The following are equivalent:
\begin{enumerate}
\item There exist a bounded linear operator 
      \(T \colon X_{L}^{\xi} \to X_{M}^{\xi}\) and 
\(\delta > 0\) such that 
\(\|T(e_{l}^{\xi})\|_0 > \delta \), for all \(l \in L\).
\item \((e_{l_n}^{\zeta})\) dominates \((e_{m_n}^{\zeta})\),
for every \(\zeta \leq \xi\).
\item \((e_{l_n}^{\xi})\) dominates \((e_{m_n}^{\xi})\).
\end{enumerate}
\end{Thm}
We recall here that a basic sequence \((x_n)\) in some Banach space \(X\)
is said to {\em dominate \/} the basic sequence 
\((y_n)\) in the Banach space \(Y\),
if there exists a constant \(C > 0\) so that
\(\|\sum_{i=1}^n a_i y_i \| \leq C
\|\sum_{i=1}^n a_i x_i \|\), for every \(n \in \mathbb{N}\)
and all scalar sequences \((a_i)_{i=1}^n\). Equivalently,
\((x_n)\) dominates \((y_n)\) if there exists a bounded linear
operator \(T\) from the closed linear span of \((x_n)\) into
the closed linear span of \((y_n)\) so that \(T(x_n)=y_n\), for
every \(n \in \mathbb{N}\). The sequences \((x_n)\) and
\((y_n)\) are {\em equivalent \/} if each one of them dominates 
the other.

As an immediate consequence of Theorem \ref{T:1} we obtain 
\begin{Cor} \label{C:1}
Let \(\xi < \omega\) and \(L=(l_n)\), \(M=(m_n)\) 
be infinite subsets of \(\mathbb{N}\). 
\begin{enumerate}
\item If \(X_L^{\xi}\) is isomorphic to a subspace of \(X_M^{\xi}\)
then \((e_{l_n}^{\xi})\) dominates 
\((e_{m_n}^{\xi})\). Consequently,
\(X_L^{\xi}\) is isomorphic to \(X_M^{\xi}\) if, and only if,
\((e_{l_n}^{\xi})\) is equivalent to \((e_{m_n}^{\xi})\). 
\item If \(X_L^{\xi}\) is isomorphic to \(X_M^{\xi}\),
then \(X_L^{\zeta}\) is isomorphic to \(X_M^{\zeta}\),
for every \(\zeta \leq \xi\).
\item Suppose that
\((e_{l_n}^{\xi})\) dominates a permutation of
\((e_{m_n}^{\xi})\). Then \((e_{l_n}^{\xi})\) dominates
\((e_{m_n}^{\xi})\).
\end{enumerate}
\end{Cor}
Theorem \ref{T:1} combined with elementary descriptive set theory
yields our next result on the structure of the subsequences of
\((e_{n}^{\xi})\), \(\xi < \omega\).
We recall here that the Banach spaces \(X\)
and \(Y\) are incomparable if neither of them is isomorphic
to a closed linear subspace of the other.
\begin{Thm} \label{T:2}
For every \(\xi < \omega\)
there exists a set \(\mathcal{A}\) (depending on \(\xi\))
consisting of infinite
subsets of \(\mathbb{N}\) and satisfying the following properties
\begin{enumerate}
\item The cardinality of \(\mathcal{A}\) is equal to the continuum.
\item For every pair \((L,M)\) of distinct elements of \(\mathcal{A}\),
the spaces \(X_{L}^\xi\) and \(X_{M}^\xi\) are incomparable.
\end{enumerate}
\end{Thm}
The proofs of the aforementioned results are given in the third
section of our paper. In the fourth section we deal with complemented
subspaces of \(X^\xi\) spanned by block bases of \((e_{n}^{\xi})\).
We show that there exists a block basis of \((e_{n}^{\xi})\) 
spanning a complemented subspace of \(X^\xi\) which 
is not isomorphic to \(X_M^\zeta\), for all
\(0 \leq \zeta \leq \xi\) and every 
infinite subset \(M\)
of \(\mathbb{N}\). We also show that there exists a block basis
of \((e_{n}^{\xi})\) spanning a subspace which is not complemented in
\(X^\xi\). 

The problem of the isomorphic classification of the complemented
subspaces of \(X^\xi\), even for block subspaces, seems rather difficult.

Part of the research for this paper was conducted while the
second author visited the University of Texas at Austin.
The second author thanks the Department of Mathematics
there, especially the Banach space group, for making the 
visit possible. Thanks are also due to Ted Odell for
several conversations regarding the results contained
herein.
\section{Preliminaries}
We shall make use of standard Banach space facts and terminology
as may be found in \cite{LT}. In this section we shall review some
of the necessary concepts. We shall also review two important hierarchies,
the Schreier hierarchy \cite{aa} and the repeated averages hierarchy
\cite{amt}. Finally we shall state some fundamental results from descriptive
set theory which will be widely used in the sequel. 
For a detailed study of descriptive set theory
we refer to \cite{ak}.

We first indicate some special notation that we will be using.
A sequence \((x_n)_{n=1}^{\infty}\) of elements of an arbitrary set  
will be conveniently denoted by \((x_n)\).
Given \(M\), a subset of \(\mathbb{N}\),   
\([M]^{< \infty}\) denotes the set of all finite subsets
of \(M\), while \([M]\) stands for the set of all infinite
subsets of \(M\). If \(M \in [\mathbb{N}]\), then the notation
\(M=(m_n)\) indicates that \(M= \{ m_1 < m_2 < \cdots \}\).
Let \(E, F\) be finite sets of integers.
We shall adopt the notation \(E < F\) to denote the relation
\(\max E < \min F\).
If \(x=\bigl(x(n)\bigr)\) belongs to \(c_{00}\), 
the space of finitely supported real valued sequences,
and \(F \in [\mathbb{N}]^{< \infty}\),
then \(x(F)= \sum_{n \in F} x(n)\), and
\(|x|(F)= \sum_{n \in F} |x(n)|\).

All Banach spaces considered throughout this paper are real.
\(\ell_1\) denotes the Banach space of the absolutely summable
sequences under the norm given by the sum of the absolute values
of the coordinates. \(c_0\) is the Banach space of the null sequences
under the norm given by the maximum of the absolute values of the
coordinates. By the term {\em ``subspace''\/} of a Banach space
we shall always mean a closed linear subspace.
A subspace \(Y\) of the Banach space \(X\)
is said to be {\em complemented \/} if it is the range of a bounded linear
projection on \(X\).
 
We next recall that if \((x_n)\) is a sequence in some normed linear space,
then the sequence \((y_n)\) is called a {\em block subsequence \/}
(resp. {\em convex block subsequence \/})
of \((x_n)\), if there exist sets \(F_i \subset \mathbb{N}\)
with \( F_1 < F_2 < \cdots \) and a sequence \((a_i)\) of scalars
(resp. non-negative scalars such that \(\sum_{n \in F_i}a_n =1\),
for every \(i \in \mathbb{N}\)) 
such that for every \(i \in \mathbb{N}\),
\(y_i = \sum_{n \in F_i}a_n x_n \).
We then denote by \(suppy_i\), the support of \(y_i\), that is the
set \(\{n \in F_i: |a_n| > 0 \}\). We shall also adopt the notation
\(y_1 < y_2 < \cdots \) to indicate that \((y_n)\) is a block
subsequence of \((x_n)\).
In case \((x_n)\) is Schauder basic, then \((y_n)\) will be called
a {\em block basis \/ } (resp. {\em convex block basis \/ })
of \((x_n)\).
 
Next we review the definition and some basic properties of the Schreier
families \(\{S_\xi\}_{\xi < \omega_1}\) \cite{aa}.
The Schreier families are defined by transfinite induction as
follows:
\[ S_0 = \bigl \{ \{n\} : n \in \mathbb{N} \} \bigr \} \cup
               \{\emptyset\}.
\]
Suppose \(S_\zeta\) has been defined for every \(\zeta < \xi\).
If \(\xi\) is a successor ordinal, say \(\xi=\zeta + 1\),
we set 
 \[  S_\xi =        \{ \cup_{i=1}^n F_i :\,  
                     n \in \mathbb{N}, \,n \leq \min F_1,  
                     F_1 < \cdots < F_n, \,     
                     F_i \in S_\zeta (i \leq n) \}
        \cup \{\emptyset\}.\]

If \(\xi\) is a limit ordinal, let \((\xi_n)\) be a preassigned
increasing sequence of successor ordinals whose limit is \(\xi\). 
We set
\[
   S_\xi = \cup_{n=1}^{\infty} \{ F \in S_{\xi_{n}}: n \leq \min F\}
         \cup \{\emptyset\}.
\]
Given \(M \in [\mathbb{N}]\) we denote by \(S_\xi[M]\) the 
family \(\{F: F \in S_\xi , \, F \subset M\}\).

An important property shared by the Schreier families
is that they are spreading: If \(\{p_1, \cdots , p_k\}
\in S_\xi\), \(p_1 < \cdots < p_k\), and \(q_1 < \cdots < q_k\)
are so that \(p_i \leq q_i\) for all \(i \leq k\), then
\(\{q_1, \cdots , q_k\} \in S_\xi\).

Of particular interest are the maximal (under inclusion)
members of \(S_\xi\). The following lemma concerning
those sets is proved in \cite{g}.
\begin{Lem} \label{L:1} 
Let \(M \in [\mathbb{N}]\) and \(\xi < \omega_1\).
Then there exists a (necessarily) unique sequence
\(\{ F_n^{\xi}(M)\}_{n=1}^{\infty}\)
of successive maximal \(S_\xi\) sets so that
\(M = \cup_{n=1}^{\infty} F_n^{\xi}(M)\).
\end{Lem}
\begin{remark}
The following stability properties of 
\(\{ F_n^{\xi}(M)\}_{n=1}^{\infty}\)
are easily verified:
\begin{enumerate}
\item If \(k_1 < k_2 < \cdots \) and 
\(N = \cup_{n=1}^{\infty} F_{k_n}^{\xi}(M)\),
then \(F_n^{\xi}(N)= F_{k_n}^{\xi}(M)\), for all
\(n \in \mathbb{N}\).
\item Let \(M=(m_i)\) and \(N=(n_i)\) be infinite
subsets of \(\mathbb{N}\). Assume that for some
\(p \in \mathbb{N}\), \(m_i = n_i\) for all
\(i \leq p\). If \(F_k^{\xi}(M)\) is contained in
\(\{m_i : i \leq p\}\), then
\(F_i^{\xi}(M) = F_i^{\xi}(N)\) for all \(i \leq k\).
\end{enumerate}
\end{remark}
In the sequel we shall make use of the following
\begin{Lem} \label{L:2}
Let \(M \in [\mathbb{N}]\), \(L \in [M]\) and \(\xi < \omega\).
Then \(\max F_1^{\xi}(M) \leq \max F_1^{\xi}(L)\).
\end{Lem}
\begin{proof}
Suppose \(L=(l_i)\) and \(M=(m_i)\). We prove the assertion
of the lemma by induction on \(\xi\). The case \(\xi=0\)
is trivial. Assume now that \(\xi \geq 1\) and that the assertion holds
for \(\xi - 1\) and all \(P\), \(Q\) with \(Q \in [P]\). 

For an arbitrary \(P \in [\mathbb{N}]\), we set \(P_1=P\) and
\[ P_i = \{ p \in P : \, p > \max F_{i-1}^{\xi - 1}(P)\} , \, i \geq 2.
\]
We observe that \(F_i^{\xi - 1}(P)=F_1^{\xi - 1}(P_i)\),
for all \(i \in \mathbb{N}\).
We also have that
\( F_1^{\xi}(P) = \cup_{i=1}^{p_1} F_1^{\xi - 1}(P_i)\),
where \(p_1 = \min P_1\).  
It follows now, by the induction hypothesis, that \(L_i \in [M_i]\),
for all \(i \in \mathbb{N}\). Therefore,
\[\max F_1^{\xi}(M)= \max F_1^{\xi - 1}(M_{m_1}) \leq
\max F_1^{\xi - 1}(L_{m_1}) \leq \max F_1^{\xi}(L)\] as \(m_1 \leq l_1\).
The proof of the lemma is now complete.
\end{proof}
We now pass to the definition of the repeated averages hierarchy
introduced in \cite{amt}. We let \((e_n)\) denote the unit vector
basis of \(c_{00}\). 
For every countable ordinal \(\xi\) and every
\(M \in [\mathbb{N}]\), we define a convex block subsequence 
\(({\xi}_{n}^{M})_{n=1}^{\infty}\)
of \((e_n)\) by transfinite induction on \(\xi\) in the following manner:
If \(\xi = 0\), then \({\xi}_{n}^{M} = e_{m_n}\), for all 
\(n \in \mathbb{N}\), where \(M = (m_n)\). 

Assume that \(({\zeta}_{n}^{M})_{n=1}^{\infty}\)
has been defined for all \(\zeta < \xi\) and \(M \in [\mathbb{N}]\).
Let \(\xi = \zeta + 1 \). Set
\[ {\xi}_{1}^{M} = \frac{1}{m_1}
                   \sum_{i=1}^{m_1} {\zeta}_{i}^{M}
\]
where \(m_1 = \min M\). Suppose that
\( {\xi}_{1}^{M} < \cdots < {\xi}_{n}^{M}\)
have been defined. Let
\[
   M_n = \{ m \in M : m > \max supp {\xi}_{n}^{M} \}
\text{ and } k_n = \min M_n .
\]
Set
\[{\xi}_{n + 1}^{M} = \frac{1}{k_n} \sum_{i=1}^{k_n}
  {\zeta}_{i}^{M_n} .
\]
If \(\xi\) is a limit ordinal, let \(({\xi}_{n} + 1)\)
be the sequence of ordinals associated to \(\xi\) in the definition
of \(S_{\xi}\), and let also
\(M \in [\mathbb{N}]\). Define 
\[
  {\xi}_{1}^{M} = [ {\xi}_{m_1} + 1]_{1}^{M}
\]
where \(m_1 = \min M\). Suppose that
\( {\xi}_{1}^{M} < \cdots < {\xi}_{n}^{M}\)
have been defined. Let
\[
   M_n = \{ m \in M : m > \max supp {\xi}_{n}^{M} \}
\text{ and } k_n = \min M_n .
\]
Set
\[
{\xi}_{n + 1}^{M} =  [ {\xi}_{k_n} + 1]_{1}^{M_n}.
\]     
The inductive definition of \(({\xi}_{n}^{M})_{n=1}^{\infty}\),
\(M \in [\mathbb{N}]\) is now complete.
The following properties are established in \cite{amt}. 

\(\mathbf{P1}\): \(({\xi}_{n}^{M})_{n=1}^{\infty}\) is a convex block
subsequence of \((e_n)\) and 
\(M = \cup_{n=1}^{\infty} supp{\xi}_{n}^{M}\)
for all \(M \in [\mathbb{N}]\) and \(\xi < {\omega}_1\).

\(\mathbf{P2}\): \(supp{\xi}_{n}^{M} \in S_{\xi}\), 
for all \(M \in [\mathbb{N}]\), \(\xi < {\omega}_1\) and \(n \in \mathbb{N}\).

\(\mathbf{P3}\): If \( M, N \in [\mathbb{N}]\), \(\xi < {\omega}_1\),
and \(supp{\xi}_{i}^{M}=supp{\xi}_{i}^{N}\),
for \(i \leq k\), then \({\xi}_{i}^{M}= {\xi}_{i}^{N}\) for \(i \leq k\).

\(\mathbf{P4}\): If \(\xi < {\omega}_1\), 
\(\{ n_k: k \in \mathbb{N}\} \subset \mathbb{N}\),
and \(\{L_k: k \in \mathbb{N}\} \subset [\mathbb{N}]\),
are such that \( supp {\xi}_{n_i}^{L_i} < supp {\xi}_{n_{i+1}}^{L_{i+1}}\),
for all \(i \in \mathbb{N}\), then letting 
\(L= \cup_{i=1}^{\infty} supp {\xi}_{n_i}^{L_i}\),
we have that \({\xi}_{i}^{L}={\xi}_{n_i}^{L_i}\), for all
\(i \in \mathbb{N}\).

Properties \(\mathbf{P3}\) and \(\mathbf{P4}\) are called stability
properties of the hierarchy 
 \(\{ ({\xi}_{n}^{M})_{n=1}^{\infty} : M \in [\mathbb{N}]\} \).
It is easily seen, by induction, that
\(supp {\xi}_n^M = F_n^{\xi}(M)\), for every \(\xi < \omega\), 
all \(M \in [\mathbb{N}]\) and \(n \in \mathbb{N}\).

In the next lemma we show that for \(\xi < \omega\) and 
\(M \in [\mathbb{N}]\) the sequence \(({\xi}_n^M)\), considered as a
sequence in \(X^\xi\), is equivalent to the unit vector basis
of \(c_0\). Moreover, the equivalence constant depends only
on \(\xi\). 
\begin{Lem} \label{L:3}
\( \|\sum_{i=1}^n {\xi}_i^M \|_{\xi} \leq \xi + 1\),
for every \(M \in [\mathbb{N}]\), \(n \in \mathbb{N}\),
and \(\xi < \omega\).
\end{Lem}
\begin{proof}
By induction on \(\xi\). The case \(\xi=0\) is trivial.
Assume the assertion holds for \(\xi - 1\). Let \(G \in S_\xi\).
We shall show that \(\sum_{i=1}^n {\xi}_i^M (G) \leq \xi + 1\),
for every \(M \in [\mathbb{N}]\) and \(n \in \mathbb{N}\).
To this end choose
\(G_1 < \cdots < G_l\), successive members of \(S_{\xi-1}\)
so that \(l \leq \min G\) and \(G=\cup_{i=1}^l G_i\).
Let also \(\{i_1, \dots i_p\}\) be an enumeration of the set
\(\{i \leq n : \, F_i^{\xi}(M) \cap G \ne \emptyset\}\).
We define 
\[ L = \cup_{t=1}^p F_{i_t}^{\xi}(M) \cup 
       \{m \in M : \, m > \max F_{i_p}^{\xi}(M)\}
\]
and observe that
\( F_{i_t}^{\xi}(M) = \cup_{j=r_{t-1} + 1}^{r_t} 
   F_j^{\xi - 1}(L)
\),
for all \(t \leq p\),
where, \(r_0 =0 < r_1 < \cdots < r_p\) are chosen so that
\(r_t - r_{t-1} = \min F_{i_t}^{\xi}(M)\), for all \(t \leq p\).
Therefore,
\[
   {\xi}_{i_t}^M = \frac{1}{\min F_{i_t}^{\xi}(M)}
   \sum_{j=r_{t-1} + 1}^{r_t} (\xi - 1)_{j}^L, \, t \leq p,
\]
and thus
\begin{align}
\sum_{t=2}^p {\xi}_{i_t}^M (G) &= \sum_{t=2}^p \sum_{s=1}^l
\frac{1}{\min F_{i_t}^{\xi}(M)}  
\sum_{j=r_{t-1} + 1}^{r_t} (\xi - 1)_{j}^L (G_s) \notag \\
&= \sum_{s=1}^l \sum_{t=2}^p    
\frac{1}{\min F_{i_t}^{\xi}(M)}  
\sum_{j=r_{t-1} + 1}^{r_t} (\xi - 1)_{j}^L (G_s) \notag \\
&\leq
\sum_{s=1}^l \frac{1}{\min F_{i_2}^{\xi}(M)} \biggl(
\sum_{t=2}^p 
\sum_{j=r_{t-1} + 1}^{r_t} (\xi - 1)_{j}^L (G_s) \biggr) \notag \\
&\leq \sum_{s=1}^l \frac{1}{\min F_{i_2}^{\xi}(M)} \xi 
\text{ (by the induction hypothesis) } \notag \\
&= \frac{l}{\min F_{i_2}^{\xi}(M)} \xi  \notag \\
&\leq \xi \text{ as } \, l \leq \min G \leq \max 
F_{i_1}^{\xi}(M) < \min F_{i_2}^{\xi}(M). \notag
\end{align}
Finally, \(\sum_{t=1}^p {\xi}_{i_t}^M (G) \leq 1 + \xi \)
and hence \(\sum_{i=1}^n {\xi}_i^M (G) \leq 1 + \xi\).
We conclude, since \(G \in S_\xi\) was arbitrary, that
\( \|\sum_{i=1}^n {\xi}_i^M \|_{\xi} \leq \xi + 1 \), as claimed.
\end{proof}
Let now \(M \in [\mathbb{N}]\).
By identifying elements of \([M]\) with
their indicator functions, \([M]\)
can be endowed with the topology of pointwise convergence. It
is not difficult to see that \([M]\) is then homeomorphic to a
\(G_{\delta}\) subset of the Cantor set \(\{0,1\}^{\mathbb{N}}\),
and thus it is a zero-dimensional Polish space. 
Further, \([M]\) is perfect
(that is contains no isolated points) and every compact subset
of \([M]\) is nowhere dense.
It is then a classical result
that \([M]\), endowed with the topology of pointwise convergence,
is homeomorphic to the space of irrational numbers with the ordinary
topology. It is worthwhile to note here that the family
\[\{ W(p_1 , \dots , p_k): k \in \mathbb{N}, \, 
   p_1 < \cdots < p_k , \, p_i \in M, \,
  i \leq k\} \] 
where \(W(p_1 , \dots , p_k)= \{L \in [M], \, L=(l_i): \,
l_i = p_i, \, i \leq k\}\), forms a basis of clopen subsets
for the topology of the pointwise convergence in \([M]\).

\section{Proofs of the main results}
This section is devoted to the proofs of Theorems \ref{T:1}
and \ref{T:2}.
\begin{Def} \label{D:1}
Let \(\xi < \omega_1\) and \(A \in [\mathbb{N}]^{< \infty}\).
We set
\[
   {\tau}_{\xi}(A) = \max \biggl\{n \in \mathbb{N}: \, A \cap
F_n^{\xi} \bigl( A \cup \{m \in \mathbb{N}: \, m > \max A \} \bigr)
 \ne \emptyset \biggr\}.
\]
\end{Def}
We observe that \({\tau}_{\xi}(A)\) remains invariant if 
\(\{m \in \mathbb{N}: \, m > \max A \}\) is replaced by
\(\{m \in M : \, m > \max A \}\), \(M \in [\mathbb{N}]\),
in Definition \ref{D:1}.
The quantity \({\tau}_{\xi}(A)\) is important for our purposes
since it will enable us state a criterion for determining whether
or not the sequence \((e_{l_n}^{\xi})\) dominates 
\((e_{m_n}^{\xi})\), where \(L=(l_n)\) and \(M=(m_n)\) belong to
\([\mathbb{N}]\). Our next lemma describes some permanence properties of
\({\tau}_{\xi}(A)\).
\begin{Lem} \label{L:21}
Let \(\xi < \omega_1\) and \(A\), \(B\) belong to 
\([\mathbb{N}]^{< \infty}\).
\begin{enumerate}
\item If \(A \subset B\) then \({\tau}_{\xi}(A) \leq 
 {\tau}_{\xi}(B)\).    
\item If \(A < B\) then 
\({\tau}_{\xi}(A \cup B) \leq {\tau}_{\xi}(A) +
  {\tau}_{\xi}(B)\).
\item If \(A=\{a_1 < \cdots < a_n\}\), 
         \(B=\{b_1 < \cdots < b_n\}\), \(n \in \mathbb{N}\),
and \(a_i \leq b_i\) for \(i \leq n\), then
\({\tau}_{\xi}(B) \leq {\tau}_{\xi}(A)\).
\item Assume that \(A=\cup_{i=1}^n A_i\),
\(B=\cup_{i=1}^n B_i\), where \(n \in \mathbb{N}\),
and \(A_1 < \cdots < A_n\), \(B_1 < \cdots < B_n\)
are maximal members of \(S_\xi\). If \(\min A_i \leq \min B_i\),
for all \(i \leq n\), then 
\({\tau}_{\xi + 1}(B) \leq {\tau}_{\xi + 1}(A)\).
\item Assume that $A = \cup_{i=1}^nA_i$ for some $n \in \mathbb{N}$. Then

\[ \tau_\xi(A) \leq (\sum^n_{i=1}\tau_\xi(A_i))(\xi+1)+1\] 

for any $\xi < \omega$.
\end{enumerate}
\end{Lem}
\begin{proof}
The first two properties are immediate consequences of
Definition \ref{D:1}. The third property follows because
\(S_\xi\) is spreading. Let us show that 4. holds.
This is accomplished by induction on \(n\). The case \(n=1\)
is easy because \({\tau}_{\xi + 1}(B) = {\tau}_{\xi + 1}(A) =1 \).
Assuming the assertion true for all \(k < n\), we set
\(k_1 = \min A_1\) and \(l_1 = \min B_1\). In case
\(l_1 \geq n\), we obtain that \(B \in S_{\xi + 1}\).
Thus \({\tau}_{\xi + 1}(B)=1\) and hence the assertion holds.   

Next suppose that \(l_1 < n\). It follows that
\(\cup_{i=1}^{k_1} A_i\) and \(\cup_{i=1}^{l_1} B_i\)
are maximal \(S_{\xi + 1}\) sets. On the other hand,
because \(n-l_1 < n\), the induction hypothesis yields
that \({\tau}_{\xi + 1}( \cup_{i=l_1 +1}^{n} B_i) \leq
{\tau}_{\xi + 1}( \cup_{i=l_1 +1}^{n} A_i)\).
But also, \(k_1 \leq l_1\) and so property 1. yields
that \({\tau}_{\xi + 1}( \cup_{i=l_1 +1}^{n} A_i) \leq
{\tau}_{\xi + 1}( \cup_{i=k_1 +1}^{n} A_i)\). The proof is
complete since \({\tau}_{\xi + 1}(B)= 1 +
{\tau}_{\xi + 1}( \cup_{i=l_1 +1}^{n} B_i)\), while
\({\tau}_{\xi + 1}(A)= 1 +
{\tau}_{\xi + 1}( \cup_{i=k_1 +1}^{n} A_i)\).

We now prove 5. Let $k = \tau_\xi(A)$ and 
$M = A \cup \{m \in \mathbb{N} : m > \max A\}$. Denote $\sum^{k-1}_{j=1}\xi^M_j$ by $x$.
By Lemma \ref{L:3}, $\|x\|_\xi \leq \xi + 1$. Hence
\begin{align*}
k - 1 &= x(A) \leq \sum^n_{i=1}x(A_i)\\
&\leq (\sum^n_{i=1}\tau_\xi(A_i))\|x\|_\xi\\
&\leq (\sum^n_{i=1}\tau_\xi(A_i))(\xi+1),
\end{align*}
from which the result follows.
\end{proof}
\begin{Def}
Let \(\xi < \omega_1\) and \(L=(l_n)\), \(M=(m_n)\) belong to
\([\mathbb{N}]\). Define
\[ d_{\xi}(L,M) = \sup \{ {\tau}_{\xi}({\phi}^{-1} A) :
\, A \in S_\xi[M] \}.\]
Where \(\phi \colon L \to M\) is the natural bijection
\(\phi(l_n)=m_n\), for all \(n \in \mathbb{N}\).
\end{Def}
The reason we introduced the quantity \(d_{\xi}(L,M)\)
is justified by our next lemma.
\begin{Lem} \label{L:22}
Let \(\xi < \omega \) and \(L=(l_n)\), \(M=(m_n)\) belong to
\([\mathbb{N}]\). Then 
\((e_{l_n}^{\xi})\) dominates 
\((e_{m_n}^{\xi})\) if and only if, 
\(d_{\xi}(L,M)\) is finite.
\end{Lem}
\begin{proof}
Suppose first that \(d_{\xi}(L,M) = p < \infty\). Let
\((a_i)_{i=1}^n\) be scalars and choose \(F \in S_\xi[M]\)
so that \(\sum_{i \in H} |a_i| =
\|\sum_{i=1}^n a_i e_{m_i}^{\xi} \|\), where we have set
\(H= \{i \leq n : \, m_i \in F \}\).
It follows, by our assumption, that we can find
\(G_1 < \cdots < G_p\) successive \(S_\xi[L]\) sets
so that \({\phi}^{-1} F \subset \cup_{j=1}^p G_j\).
We now set \(H_j = \{ i \in H : \, l_i \in G_j \}\),
for all \(j \leq p\). It is clear that 
\(H = \cup_{j=1}^p H_j\) and moreover, 
\(\{l_i : \, i \in H_j\}\) belongs to \(S_\xi[L]\).
Finally, 
\[ \sum_{i \in H} |a_i| = \sum_{j=1}^p \sum_{i \in H_j}
   |a_i| \leq p \, \biggl\|\sum_{i=1}^n a_i e_{l_i}^{\xi} \biggr\|.
\]
Thus, \(\|\sum_{i=1}^n a_i e_{m_i}^{\xi} \| \leq p \,  
      \|\sum_{i=1}^n a_i e_{l_i}^{\xi} \|\).

Conversely, assume that 
\((e_{l_n}^{\xi})\) \(C\)-dominates 
\((e_{m_n}^{\xi})\) and let \(F \in S_\xi[M]\).
Suppose that \({\tau}_{\xi}({\phi}^{-1} F)=k\). 
It follows that there exist \(G_1 < \cdots < G_{k-1}\),
successive maximal \(S_\xi[L]\) sets so that
\( \cup_{i=1}^{k-1} G_i \subset {\phi}^{-1} F\).
Put \[Q= \cup_{i=1}^{k-1} G_i \cup 
\{ l \in L : \, l > \max G_{k-1}\}. \]
We may write \({\xi}_i^Q = \sum_{j \in G_i} 
a_j^i e_j^{\xi}\), with \(\sum_{j \in G_i} a_j^i =1\),
for all \(i \leq k-1\).
If we apply Lemma \ref{L:3}, we obtain
\begin{align}
C(\xi + 1) &\geq C \, \biggl \|\sum_{i=1}^{k-1} {\xi}_i^Q 
\biggl\|_{\xi}
= C \, \biggl\|\sum_{i=1}^{k-1}  
\sum_{j \in G_i} a_j^i e_j^{\xi} \biggr\| \notag \\
&\geq  \biggl\|\sum_{i=1}^{k-1}  
\sum_{j \in G_i} a_j^i e_{\phi(j)}^{\xi} \biggr\|
\geq k-1 , \notag
\end{align}
as \(\cup_{i=1}^{k-1} \{ \phi(j) : j \in G_i\} 
\subset F\) and \(\sum_{j \in G_i} a_j^i =1\).
Hence, \(k \leq C(\xi + 1) + 1\) which implies that
\(d_{\xi}(L,M) \leq C(\xi + 1) + 1\) as \(F\) was
an arbitrary \(S_\xi[M]\) set.
\end{proof}
We shall next show that \((e_{n}^{\xi})\) 
has ``many'' non-equivalent subsequences.
\begin{Lem} \label{L:23}
Let \(1 \leq \xi < \omega\), \(N \in [\mathbb{N}]\) and set
\[\mathcal{D}=\{ (L,M) \in [N] \times [N] : \,
d_{\xi}(L,M)=d_{\xi}(M,L)=\infty\}.\]
Then \(\mathcal{D}\) is a \(G_{\delta}\) dense subset
of \([N] \times [N]\).
\end{Lem}
\begin{proof}
By Baire's theorem, it suffices to show that the sets 
\(\{(L,M) \in [N] \times [N] : \, d_{\xi}(L,M) < \infty\}\)
and 
\(\{(L,M) \in [N] \times [N] : \, d_{\xi}(M,L) < \infty\}\) 
are first category \(F_{\sigma}\) subsets of
\([N] \times [N]\). Indeed, we may write
\[
  \{(L,M) \in [N] \times [N] : \, d_{\xi}(L,M) < \infty\}=
 \cup_{n=1}^{\infty} 
\{(L,M): \, d_{\xi}(L,M) \leq n\}.
\]
It is easy to see that each set in the union is closed
in \([N] \times [N]\) and thus it remains to show that
\(\{(L,M) \in [N] \times [N] : \, d_{\xi}(L,M) \leq n\}\)
has empty interior in \([N] \times [N]\).
If that were not the case, choose \(\mathcal{U}\) and
\(\mathcal{V}\), non-empty basic clopen subsets of
\([N]\) so that \(\mathcal{U} \times \mathcal{V}\)
is contained in 
\(\{(L,M) \in [N] \times [N] : \, d_{\xi}(L,M) \leq n\}\).
There exist \(p_1 < \dots < p_k\) in \(N\) so that
\(\mathcal{V}= W(p_1, \dots , p_k)\). Fix \(L \in \mathcal{U}\),
\(L=(l_i)\). If \(M \in [N]\), \(\min M > p_k\), let
\(P=\{p_1, \cdots , p_k\} \cup M\). Since 
\(d_\xi(L,P) \leq n\),
it follows that if \(L_k = (l_{k+i})_{i=1}^{\infty}\), then
\(d_{\xi}(L_k,M) \leq n\).
By Lemma \ref{L:22} this 
implies that \((e_l^{\xi})_{l \in L_k}\) is equivalent to
the unit vector basis of \(\ell_1\) which is absurd.
Arguing similarly, we also obtain that 
\(\{(L,M) \in [N] \times [N] : \, d_{\xi}(M,L) < \infty\}\)
is first category, 
\(F_{\sigma}\) subset of
\([N] \times [N]\). 
\end{proof}
We also need the following result which is
a special case of a theorem by Mycielski \cite{m} 
and Kuratowski \cite{k}
(cf. also \cite{ak}).
\begin{Prop} \label{P:1}
Let \(K\) be a perfect Polish space and \(G\) a \(G_{\delta}\)
dense subset of \(K \times K\). There exists \(C\), a subset
of \(K\) homeomorphic to the Cantor set such that
\(C \times C \setminus \Delta \subset G\) (here
\(\Delta\) is the diagonal subset of \(K \times K\)).
\end{Prop}
This result may be found in \cite{ak} (p. 129, Theorem 19.1)
but we shall include
a proof to be thorough.
\begin{Lem} \label{L:24}
Let \(K\) be Polish and \(G\) be an open dense subset of \(K \times K\).
Let also \((A_i)_{i=1}^n\) (\(n \geq 2\)) be a finite sequence
of open non-empty subsets of \(K\). Then for every
\(\epsilon > 0\) there exist
\((B_i)_{i=1}^n\), open non-empty subsets of \(K\),
so that \(diam B_i < \epsilon\), for all \(i \leq n\) and
\(\overline{ B_i } \times \overline{ B_j }
\subset G \cap (A_i \times A_j)\),
for all \(i \ne j\) in \(\{1, \dots , n\}\).
\end{Lem}
\begin{proof}
By induction on \(n\). Suppose first that \(n=2\).
Since \((A_1 \times A_2) \cap G \ne \emptyset\),
there exist \(C_1\), \(C_2\), open non-empty subsets
of \(K\) whose diameters are smaller than \(\epsilon\),
so that \(\overline{ C_1 } \times \overline{ C_2 }
\subset G \cap (A_1 \times A_2)\).
Further, \((C_2 \times C_1) \cap G \ne \emptyset\)
and thus there exist \(B_1\), \(B_2\), open non-empty
subsets of \(K\), so that 
\(\overline{ B_2 } \times \overline{ B_1 }
\subset G \cap (C_2 \times C_1)\). Of course 
\(B_1\) and \(B_2\) satisfy the conclusion
of the lemma for \(n=2\).

Next assume \(n > 2\) and that the result holds for
\(n-1\).  We can therefore choose
\((C_i)_{i=1}^{n-1}\), open non-empty subsets of
\(K\) with diameters smaller than \(\epsilon\),
so that 
\(\overline{ C_i } \times \overline{ C_j }
\subset G \cap (A_i \times A_j)\), 
for all \(i \ne j\) in \(\{1, \dots , n-1\}\).
Next, set \(A_{n,0}=A_n\) and choose, as in the case
\(n=2\), \((B_i)_{i=1}^{n-1}\),
\((A_{n,i})_{i=1}^{n-1}\), non-empty open subsets of \(K\)
with diameters smaller than \(\epsilon\),
so that 
\(\overline{ B_i } \times \overline{ A_{n,i} }
\subset G \cap (C_i \times A_{n,i-1})\) and
\(\overline{ A_{n,i} } \times \overline{ B_i }
\subset G \cap (A_{n,i-1} \times C_i )\),
for all \(i \leq n-1\).
Set \(B_n = A_{n,n-1}\) and it is easy to check that
\((B_i)_{i=1}^n\) is the desired sequence. 
\end{proof}
\begin{proof}[Proof of Proposition \ref{P:1}]
Since \(K\) contains no isolated points,
\(\Delta\) is nowhere dense in \(K \times K\).
Hence, \(G \cap (K \times K \setminus \Delta)\)
is a \(G_{\delta}\) dense subset of \(K \times K \).
We shall therefore assume, without loss of generality,
that \(G \cap \Delta = \emptyset\). Now let 
\((G_n)\) be a decreasing sequence of open dense subsets
of \(K \times K\), whose intersection is \(G\). We can
assume that \(G_n \cap \Delta = \emptyset\), for all
\(n \in \mathbb{N}\).

We shall construct a collection
\(\bigl\{U_{\alpha}: \, \alpha \in \{0,1\}^n, \, 
n \in \mathbb{N} \bigr\}\)
of open non-empty subsets of \(K\) so that the following
properties are satisfied for every \(n \in \mathbb{N}\):

(i) \(\overline{ U_\alpha } \cap \overline{ U_\beta }
= \emptyset\), whenever \(\alpha \ne \beta\) in 
\(\{0,1\}^n\).

(ii) \(\overline{ U_\alpha } \subset U_\beta\),
for all \(\alpha \in \{0,1\}^n\) and every
\(\beta \in \{0,1\}^m\), (\(m < n\)), initial segment
of \(\alpha\). 

(iii) \(diam U_\alpha < \frac{1}{n}\), for every 
\(\alpha \in \{0,1\}^n\).

(iv) \(\overline{ U_\alpha } \times 
      \overline{ U_\beta } \subset G_n \),
whenever \(\alpha \ne \beta\) in \(\{0,1\}^n\).

Once this is accomplished, we let
\[ C = \bigl\{ x \in K : \, \exists \, \alpha \in 
\{0,1\}^{\mathbb{N}},\,   \{x\}= 
\cap_{n=1}^{\infty} U_{\alpha | n} \bigr \}.
\]
Where \(\alpha | n = (a_1, \dots , a_n)\),
if \(\alpha=(a_i) \in \{0,1\}^{\mathbb{N}}\).
It is a standard result that \(C\) is homeomorphic to
the Cantor set. Property (iv) yields that \(C\)
satisfies the conclusion of Proposition \ref{P:1}.

The construction is done by induction on 
\(n \in \mathbb{N}\). For \(n=1\) choose
\(W_0\) and \(W_1\), open non-empty subsets of
\(K\) so that \(W_0 \times W_1 \subset G_1\).
\(W_0\) and \(W_1\) are disjoint since
\(G_1 \cap \Delta = \emptyset\). If we apply
Lemma \ref{L:24}, for \(\epsilon=1\),
on the dense open subset \(G_1\) and the open
sets \(W_0\) and \(W_1\), we shall find
\(U_0\), \(U_1\),
non-empty open subsets of \(K\), satisfying
properties (i)-(iv), for \(n=1\).

Now suppose that for every \(k \leq n\) we have constructed
\(\bigl\{U_{\alpha}:\, \alpha \in \{0,1\}^k \bigr\}\),
a collection 
of open non-empty subsets of \(K\)
whose members satisfy properties (i)-(iv),
for \(k\). 
Let \(\{d_1, \dots , d_p\}\), \(p=2^n\),
be an enumeration of \(\{0,1\}^n\).
Another application of Lemma \ref{L:24}
yields \(W_{j0}\) and \(W_{j1}\), 
non-empty open subsets of \(K\), \(j \leq p\), so that
\( \overline{ W_{jr} } \times 
   \overline{ W_{js} } \subset
  ( U_{d_j} \times U_{d_j}) \cap G_{n+1}\),
for every \(j \leq p\) and all pairs \((r,s)\)
of distinct elements of \(\{0,1\}\).
It follows, since \(G_{n+1} \cap \Delta = \emptyset\),
that \(W_{j0} \cap W_{j1} = \emptyset\).
According to the induction hypothesis
\( \overline{ U_{d_j} } \cap 
           \overline{ U_{d_i} } = \emptyset\),
for all \(i \ne j\) in \(\{1, \dots , p\}\),
and thus 
\( \overline{ W_{jr} } \cap
   \overline{ W_{is} } = \emptyset\),
for all \((j,r) \ne (i,s)\) in
\(\{1, \dots , p\} \times \{0,1\}\).

We next apply Lemma \ref{L:24}, for \(
\epsilon = \frac{1}{n+1}\),
on the family \(\bigl \{ W_{jr}: \, (j,r) \in  
\{1, \dots , p\} \times \{0,1\} \bigr \}\)
and the dense open subset \(G_{n+1}\).
We shall obtain 
\(\bigl \{ U_{\alpha}: \, \alpha \in \{0,1\}^{n+1}
 \bigr \}\), a collection of non-empty open subsets
of \(K\), so that 
\(\overline{ U_{\alpha} } \times 
  \overline{ U_{\beta} } 
\subset (W_{jr} \times W_{is}) \cap G_{n+1} \)
whenever \(\alpha = (d_j , r)\), 
\(\beta = (d_i , s)\)
and \((j,r) \ne (i,s)\) in
\(\{1, \dots , p\} \times \{0,1\}\).
Evidently, \(\bigl \{ U_{\alpha}: \, \alpha \in \{0,1\}^{n+1}
 \bigr \}\) satisfies properties (i)-(iv). 
The inductive step as well as the proof of
the proposition are now complete.
\end{proof}

Assuming we have proved Theorem \ref{T:1} 
and Corollary \ref{C:1}, let us now show 
how to derive Theorem \ref{T:2} from our previously
obtained results.
\begin{proof}[Proof of Theorem \ref{T:2}]
Let \(\mathcal{D}\) be as in the statement of
Lemma \ref{L:23}, where we have taken \(N=\mathbb{N}\).
We can apply Proposition \ref{P:1} for the space
\([\mathbb{N}]\) and the subset \(\mathcal{D}\)
to obtain \(\mathcal{A} \subset [\mathbb{N}]\), homeomorphic
to the Cantor set and such that
\(\mathcal{A} \times \mathcal{A} \setminus \Delta
\subset \mathcal{D}\). Lemma \ref{L:22} and Corollary
\ref{C:1} yield that \(\mathcal{A}\)
is the desired subset of \([\mathbb{N}]\).
\end{proof} 
We shall next pass to the proof of Theorem \ref{T:1}.
We first prove some necessary lemmas.
\begin{Lem} \label{L:25}
Let \(G \in [\mathbb{N}]^{ < \infty}\) 
and \(\xi < \omega\). The following are equivalent:
\begin{enumerate}
\item \(G\) is a member (resp. maximal member) of \(S_\xi\).
\item For every \(0 \leq \zeta \leq \xi\) there exist
\(n \in \mathbb{N}\) and \(G_1 < \cdots < G_n\) successive
members (resp. maximal members) of \(S_\zeta\) so that
\(G= \cup_{i=1}^n G_i\) and 
\(\{\min G_i: \, i \leq n\}\) is a member (resp. maximal member)
of \(S_{\xi - \zeta}\).
\item There exist \(0 \leq \zeta \leq \xi\),
\(n \in \mathbb{N}\) and \(G_1 < \cdots < G_n\) successive
members (resp. maximal members) of \(S_\zeta\) so that
\(G= \cup_{i=1}^n G_i\) and 
\(\{\min G_i: \, i \leq n\}\) is a member (resp. maximal member)
of \(S_{\xi - \zeta}\).
\end{enumerate}
\end{Lem}
\begin{proof}
We show that all three conditions are equivalent
for the members of \(S_\xi\).

\(1. \Rightarrow 2.\) By induction on \(\xi\).
If \(\xi = 0\) the assertion is trivial. Suppose now
\(\xi \geq 1\) and that the assertion holds for
\(\xi -1\). Let \(\zeta \leq \xi\). If \(\zeta = \xi\),
the assertion is again trivial. So assume \(\zeta < \xi\).
Choose \(H_1 < \cdots < H_p\) in \(S_{\xi - 1}\) so that
\(p \leq \min H_1\) and \(G = \cup_{i=1}^p H_i\). The induction
hypothesis yields that for each \(i \leq p\) there exist
\(H_{i1} < \cdots < H_{ir_i}\) in \(S_\zeta\) so that
\(\{\min H_{ij}: \, j \leq r_i \}\) belongs to
\(S_{\xi - \zeta -1}\) and
\(H_i = \cup_{j=1}^{r_i} H_{ij}\).
Let \(\{G_1, \dots , G_n\}\) be an enumeration of the
set \(\{H_{ij}: \, j \leq r_i , \, i \leq p\}\)
so that \(G_1 < \cdots < G_n\). Note that
\(\{\min G_i: \, i \leq n\}= \cup_{i=1}^p 
\{\min H_{ij}: \, j \leq r_i\}\) 
and so it is a member 
of \(S_{\xi - \zeta}\) as \(p \leq \min H_{11}=\min H_1\).

\(2. \Rightarrow 3.\) This implication is trivial.

\(3. \Rightarrow 1.\)  
By induction on \(\xi\).
If \(\xi = 0\) the assertion is trivial. Suppose now
\(\xi \geq 1\) and that the assertion holds for
\(\xi -1\). Let \(\zeta \leq \xi\). If \(\zeta = \xi\),
the assertion is again trivial. So assume \(\zeta < \xi\).
We first apply \(1. \Rightarrow 2.\) for the set
\(\{\min G_i: \, i \leq n\} \in S_{\xi - \zeta}\)
to obtain \(H_1 < \cdots < H_p\) in
\( S_{\xi - \zeta - 1}\) so that 
\(\{\min G_i: \, i \leq n\}= \cup_{i=1}^p H_i\)
and \(\{\min H_i: \, i \leq p\} \in S_1\).
Set \(I_i = \{j \leq n: \min G_j \in H_i\}\), \(i \leq p\).
Then \(\cup_{j \in I_i} G_j \in S_{\xi -1}\),
by the induction hypothesis since 
\(\{\min G_j: \, j \in I_i\}=H_i\) which belongs to
\( S_{\xi - \zeta - 1}\).
Finally, \(G= \cup_{i=1}^p (\cup_{j \in I_i} G_j )
\in S_\xi\), as \(\min G = \min H_1 \geq p\).
The latter inequality holds because 
\(\{\min H_i: \, i \leq p\} \in S_1\).

The proof for the case of maximal Schreier sets
requires only minor modifications. Namely, all the
sets which belong to an appropriate class 
\(S_\alpha\), \( \alpha \leq \xi\)
and appear in the previous arguments,
can be taken to be maximal members of \(S_\alpha\).
\end{proof}
\begin{Lem} \label{L:26}
Let \(\xi < \omega_1\). Suppose that \(L=(l_i)\),
\(M=(m_i)\) belong to \([\mathbb{N}]\) and satisfy
\(l_i < m_i < l_{i+1}\), for every \(i \in \mathbb{N}\).
Then \((e_{l_i}^{\xi})\) is \(2\)-equivalent to
\((e_{m_i}^{\xi})\). That is,
\(\|\sum_{i=1}^n a_i e_{l_i}^{\xi}\| \leq 
  \|\sum_{i=1}^n a_i e_{m_i}^{\xi}\| \leq 2
  \|\sum_{i=1}^n a_i e_{l_i}^{\xi}\|\), 
for every \(n \in \mathbb{N}\) and all scalar sequences
\((a_i)_{i=1}^n\).
\end{Lem}
We omit the easy proof and pass to
\begin{Lem} \label{L:27}
Let \(\xi < \omega\) and \(0 \leq \zeta \leq \xi\).
Then for every \(L \in [\mathbb{N}]\),
\(({\zeta}_n^L)\), considered as a sequence in
\(X^\xi\), is \(12(\zeta + 1)\)-equivalent to
\((e_{q_n}^{\xi - \zeta})\). Here we have set
\(q_n = \min F_n^{\zeta}(L)\), for all \(n \in \mathbb{N}\).
\end{Lem}
\begin{proof}
Let \(n \in \mathbb{N}\) and \((a_i)_{i=1}^n\)
be scalars. Choose 
\(G \subset \{q_1, \dots , q_n\}\) with
\(G \in S_{\xi - \zeta}\) such that
\(\sum_{i \in I} |a_i|= \|\sum_{i=1}^n a_i 
   e_{q_i}^{\xi - \zeta}\|\), where we have set
\(I=\{ i \leq n: q_i \in G\}\). It follows that
\(H=\cup_{i \in I} F_i^{\zeta}(L) \in S_\xi\), by Lemma
\ref{L:25}. Hence,
\(\|\sum_{i=1}^n a_i {\zeta}_i^L\|_{\xi} 
\geq \sum_{i \in I} |a_i| \) and thus
\(\|\sum_{i=1}^n a_i 
   e_{q_i}^{\xi - \zeta}\| \leq 
   \|\sum_{i=1}^n a_i {\zeta}_i^L\|_{\xi}\).

We next show that 
\(\|\sum_{i=1}^n a_i {\zeta}_i^L\|_{\xi} \leq
12(\zeta + 1) \|\sum_{i=1}^n a_i 
   e_{q_i}^{\xi - \zeta}\|\).
Let \(G \in S_\xi[L]\) be maximal and put
\(G_i = F_i^{\zeta}(L) \cap G\), for all
\(i \leq n\).
We apply Lemma \ref{L:25} to find \(p \in \mathbb{N}\)
and \(H_1 < \cdots < H_p\) maximal members of
\(S_\zeta\) so that \(G = \cup_{j=1}^p H_j\)
and \(\{\min H_j: \, j \leq p\}\) is a maximal member
of \(S_{\xi - \zeta}\).

We claim that each of the \(G_i\)'s can intersect at
most two of the \(H_j\)'s. Indeed, assume that for
some \(i\) and \(j_1 < j_2 < j_3\) we had that
\(G_i \cap H_{j_r} \ne \emptyset\), for all
\(r \leq 3\). Then \(H_{j_2} \subset G_i\)
because \(H_{j_2} \subset [\min G_i , \max G_i]\).
Thus, \(H_{j_2} \subset F_i^{\zeta}(L)\) and hence
\(H_{j_2}=F_i^{\zeta}(L)\), by the maximality of
\(H_{j_2}\). It follows that \(H_{j_2}=G_i\)
which is a contradiction as \(H_{j_1} \cap H_{j_2} 
=\emptyset\).

Therefore our claim holds and evidently, for each
\(i \leq n\), \(G_i\) intersects either exactly one
of the \(H_j\)'s, or exactly two (consecutive)
\(H_j\)'s. We can thus partition \(\{1, \dots , n\}\)
in the following two subsets:
\[ I_1 = \{ i \leq n: \, G_i \subset H_j, \, \text{
            for some } j \leq p\},\]
\[ I_2 = \{ i \leq n: \, \exists j_1 < j_2 \leq p,\, 
G_i \subset H_{j_1} \cup H_{j_2}, \, 
G_i \cap H_{j_r} \ne \emptyset, \, r \leq 2\}.\]
Let 
\(k_i = \max F_i^{\zeta}(L)\), for all \(i \leq n\).
We now have the following 

{\em CLAIM}. Suppose that 
\(T_i \subset F_i^{\zeta}(L)\),
for all \(i \leq n\). Assume also that for each
\(i \leq n\) there exists \(j \leq p\) so that
\(T_i \subset H_j\). Then
\[\biggl|\sum_{i=1}^n a_i {\zeta}_i^L (\cup_{m=1}^n T_m) \biggr|
   \leq (\zeta + 1) \biggl\|\sum_{i=1}^n a_i 
   e_{k_i}^{\xi - \zeta} \biggr\|.
\]
Once the claim is established we finish
the proof as follows: Observe that our claim yields
\[\biggl|\sum_{i \in I_1} a_i {\zeta}_i^L (G) \biggr|
   \leq (\zeta + 1) \biggl\|\sum_{i=1}^n a_i 
   e_{k_i}^{\xi - \zeta} \biggr\|.\]
On the other hand, if \(i \in I_2\) there exist
\(A_i < B_i\) so that \(G_i = A_i \cup B_i\)
and each element of \(\{A_i , B_i\}\) is contained in
some \(H_j\). Our claim then yields that
\[\biggl|\sum_{i \in I_2} a_i {\zeta}_i^L (G) \biggr|
   \leq 2(\zeta + 1) \biggl\|\sum_{i=1}^n a_i 
   e_{k_i}^{\xi - \zeta} \biggr\|.\]
Therefore, \( |\sum_{i=1}^n a_i {\zeta}_i^L (G)|
   \leq 3(\zeta + 1) \|\sum_{i=1}^n a_i 
   e_{k_i}^{\xi - \zeta}\|\).
It follows now, since \(G\) was arbitrary, that
\(\|\sum_{i=1}^n a_i {\zeta}_i^L\|_{\xi} \leq
 6(\zeta + 1) \|\sum_{i=1}^n a_i 
   e_{k_i}^{\xi - \zeta}\|\).
The desired estimate follows now from Lemma \ref{L:26}.

We proceed now to prove our claim.
Let \(R_j = \{i \leq n: T_i \ne \emptyset, \, 
              T_i \subset H_j\}\), \(j \leq p\),
and choose \(i_j \in R_j\) such that
\(\max_{r \in R_j} |a_r| = |a_{i_j}|\).
\begin{align}
\biggl|\sum_{i=1}^n a_i {\zeta}_i^L (\cup_{m=1}^n T_m) \biggr| 
&=
\biggl|\sum_{i=1}^n a_i {\zeta}_i^L (T_i) \biggr| =
\biggl|\sum_{j=1}^p \sum_{i \in R_j} a_i {\zeta}_i^L (T_i) \biggr|
\notag \\
&\leq \sum_{j=1}^p \biggl|\biggl(
 \sum_{i \in R_j} a_i {\zeta}_i^L \biggr)
 \bigl(\cup_{m \in R_j} T_m \bigr) \biggr| \leq
 \sum_{j=1}^p \biggl\|\sum_{i \in R_j} a_i {\zeta}_i^L 
 \biggr\|_{\zeta}, \notag \\
&\text{ because } \cup_{m \in R_j} T_m \in S_{\zeta},
\notag \\
&\leq \sum_{j=1}^p (\zeta + 1) \max_{i \in R_j} |a_i|,
 \text{ by Lemma \ref{L:3} }, \notag \\
&= \sum_{j=1}^p (\zeta + 1) |a_{i_j}| \leq
   (\zeta + 1) \biggl\|\sum_{i=1}^n a_i 
   e_{k_i}^{\xi - \zeta} \biggr\|. \notag
\end{align}
The last inequality holds since
\(T_{i_j} \subset H_j\) implies that 
\(\min H_j \leq k_{i_j}\), for all \(j \leq p\)
and thus \(\{k_{i_j}: \, j \leq p\} \in
S_{\xi - \zeta}\).
The proof of the lemma is now complete.
\end{proof}
We recall here that a sequence \((x_n)\) in some
Banach space is said to be an \({\ell}_1^\xi\)-{\em
spreading model}, \(\xi < \omega_1\), provided that
there exists a constant \(C > 0\) so that
\(\|\sum_{i \in F} a_i x_i \| \geq C 
\sum_{i \in F} |a_i|\), for every \(F \in S_\xi\)
and all scalars \((a_i)_{i \in F}\).
\begin{remark}
It is easy to see that every subsequence of
\((e_n^{\xi})\) is an \({\ell}_1^\xi\)-spreading model
in \(X^\xi\). However, Lemma \ref{L:3} implies that
no subsequence of \((e_n^{\xi})\) is an
\(\ell_1^{\xi + 1}\)-spreading model in \(X^\xi\).
\end{remark}
\begin{Prop} \label{P:2}
Suppose \(L=(l_i)\), \(M=(m_i)\) belong to 
\([\mathbb{N}]\) and that \(\xi < \omega\).
Assume further that there exist a map
\(\psi \colon L \to M\) and a bounded  
linear operator \(T \colon X_L^\xi \to
X_M^\xi \) so that \(T(e_l^\xi) = e_{\psi(l)}^\xi\),
for every \(l \in L\). Then there exists an integer \(D > 0\)
so that \({\tau}_{\zeta}({\psi}^{-1} F) \leq D\),
for every \(F \in S_\zeta[M]\) and all
\(0 \leq \zeta \leq \xi\).
\end{Prop}
\begin{proof}
We first note that \({\psi}^{-1} F \in [L]^{< \infty}\),
for every \(F \in [M]^{< \infty}\). Indeed, if that were
not the case, we would find \(m \in M\) and \(N \in [L]\)
so that \(\psi(l)=m\), for every \(l \in N\). It follows
that \(T({\xi}_i^N)= e_m^\xi\), for all \(i \in \mathbb{N}\).
But this contradicts Lemma \ref{L:3} because \(T\) is
bounded.

Fix \(0 \leq \zeta \leq \xi\). Our first task is to show that
\(\sup_n {\tau}_{\zeta} \bigl({\psi}^{-1} F_n^{\zeta}(P) \bigr)
< \infty\), for every \(P \in [M]\).
Suppose this is not the case and so
\(\sup_n {\tau}_{\zeta} \bigl({\psi}^{-1} F_n^{\zeta}(P) \bigr)
 =\infty\), for some \(P \in [M]\). We claim that
there exist a sequence of positive integers, \((n_i)\),
and a sequence
of successive maximal \(S_{\zeta + 1}[L]\)-sets,
\((G_i)\), so that letting \(q_i = \min G_i\),
for all \(i \in \mathbb{N}\), the following is satisfied:
\[ G_i \setminus \{q_i\} \subset
   {\psi}^{-1} F_{n_i}^{\zeta}(P), \text{ for all }
   i \in \mathbb{N}.\]
Indeed, choose \(n_1\) so that 
\({\tau}_{\zeta} \bigl({\psi}^{-1} F_{n_1}^{\zeta}(P) \bigr)
> l_1\). Put \(q_1 = l_1\). Because 
\({\psi}^{-1} F_{n_1}^{\zeta}(P)\) contains at least \(l_1\)
successive maximal \(S_\zeta[L]\)-sets, it is clear that there
exists \(H_1 \subset {\psi}^{-1} F_{n_1}^{\zeta}(P)\),
\(q_1 < \min H_1\),
so that \(G_1 = \{q_1\} \cup H_1 \) is
a maximal \(S_{\zeta + 1}[L]\)-set.

Put \(l_{t_1}= \max G_1\) and 
\(w_1 = {\tau}_{\zeta}(\{l_1, \dots , l_{t_1}\})\).
We can find \(n_2 > n_1\) so that 
\({\tau}_{\zeta} \bigl({\psi}^{-1} F_{n_2}^{\zeta}(P) \bigr) 
> w_1 + l_{t_1 + 1}\). Now,
\(\{l \in L: \, l \geq l_{t_1 + 1}\} \cap 
   {\psi}^{-1} F_{n_2}^{\zeta}(P)\)
must contain at least \(l_{t_1 + 1}\) successive maximal
\(S_\zeta[L]\)-sets. If not, then
\({\tau}_{\zeta} [\{l \in L: \, l \geq l_{t_1 + 1}\} \cap 
   {\psi}^{-1} F_{n_2}^{\zeta}(P)] \leq l_{t_1 + 1}\)
and thus \({\tau}_{\zeta}({\psi}^{-1} F_{n_2}^{\zeta}(P))
\leq l_{t_1 + 1} + w_1\), by Lemma \ref{L:21}. But this
contradicts the choice of \(n_2\).

We set \(q_2 = l_{t_1 + 1}\) and arguing as we did in 
the case \(i=1\), we can find 
\(H_2 \subset {\psi}^{-1} F_{n_2}^{\zeta}(P)\),
\(q_2 < \min H_2\),
so that \(G_2 = \{q_2\} \cup H_2\) is a maximal
\(S_{\zeta + 1}[L]\)-set.
We next put \(l_{t_2} = \max G_2\) and continue in the same
fashion to obtain sequences \((n_i)\), \((G_i)\)
satisfying the desired properties.

Let \(Q= \cup_{i=1}^{\infty} G_i\). Clearly,
\(Q \in [L]\) and \(F_i^{\zeta + 1}(Q)=G_i\),
for all \(i \in \mathbb{N}\).
We now set \(R= \cup_{i=1}^{\infty} F_{n_i}^{\zeta}(P)\).
Then, \(R \in [M]\) and 
\(F_i^{\zeta}(R)= F_{n_i}^{\zeta}(P)\),
for all \(i \in \mathbb{N}\).
We observe that if \(q \in G_i \setminus \{q_i\}\),
then \(T(e_q^\xi)= e_m^\xi\), for some 
\(m \in F_i^{\zeta}(R)\).

Next write \((\zeta + 1)_i^Q = a_i e_{q_i}^\xi + (1-a_i)u_i\),
for all \(i \in \mathbb{N}\). Here, \(u_i\) is a convex combination
of the vectors \((e_q^\xi)_{q \in G_i \setminus \{q_i\}}\)
and \(0 < a_i < 1\). Evidently, \(\lim_i a_i = 0\).
Observe that \(Tu_i\) is a convex combination
of the vectors \((e_m^\xi)_{m \in F_i^\zeta(R)}\)
and thus \(\|Tu_i\|_\xi =1\), for all \(i \in \mathbb{N}\).

It must be the case that \(\zeta < \xi\) for if not,
Lemma \ref{L:3} yields \(\lim_i \|u_i\|_\xi = 0\).
On the other hand 
\(\|Tu_i\|_\xi =1\), for all \(i \in \mathbb{N}\).
Hence \(T\) is not bounded contrary to our assumption.
Therefore, \(\zeta < \xi\) and so \(\|u_i\|_\xi = 1\),
for all \(i \in \mathbb{N}\).

Recall that \( \lim_i \|(1-a_i)Tu_i\|_\xi = 1\) and 
\((1-a_i)Tu_i\) is supported by \(F_i^\zeta(R)\).
Using Lemma \ref{L:25}, it is easy to check that 
\(((1-a_i)Tu_i)\) is an \(\ell_1^{\xi - \zeta}\)-spreading
model in \(X_M^\xi\), and consequently, since
\(\lim_i a_i = 0\), \((T[(\zeta + 1)_i^Q])\) is also
an \(\ell_1^{\xi - \zeta}\)-spreading
model in \(X_M^\xi\). We conclude, as \(T\) is bounded,
that \(((\zeta + 1)_i^Q)\) is an 
\(\ell_1^{\xi - \zeta}\)-spreading
model in \(X^\xi\). However, if we apply Lemma \ref{L:27}
we obtain that \((e_{q_n}^{\xi - \zeta -1 })\) is
an \(\ell_1^{\xi - \zeta}\)-spreading
model in \(X^{\xi - \zeta - 1}\). But this contradicts
with the remark after Lemma \ref{L:27}.
Hence,
\(\sup_n {\tau}_{\zeta} \bigl({\psi}^{-1} F_n^{\zeta}(P) \bigr)
< \infty\), for every \(P \in [M]\). 
It follows that
\[[M] = \cup_{k=1}^{\infty} \{ P \in [M]: \, 
         {\tau}_{\zeta} \bigl({\psi}^{-1} F_n^{\zeta}(P)\bigr)
         \leq k, \, \forall n \in \mathbb{N}\}.\]
It is easily seen that every set in the union is closed in
\([M]\). Baire's theorem now yields \(k_\zeta \in \mathbb{N}\)
and \(r_1^\zeta < \cdots < r_{s_\zeta}^\zeta\) 
in \(M\) so that if
\(P \in [M]\), \(P=(p_i)\), and \(p_i=r_i^\zeta\), \(i \leq s_\zeta\),
then \( {\tau}_{\zeta} \bigl({\psi}^{-1} F_n^{\zeta}(P)\bigr)
         \leq k_\zeta\), for all \(n \in \mathbb{N}\).
It follows now that there exists \(m_0^\zeta \in M\) so that if
\(F \in S_\zeta[M]\), \(\min F > m_0^\zeta\), then
\({\tau}_{\zeta}({\psi}^{-1} F) \leq k_\zeta\).

Finally, choose \(D_\zeta \in \mathbb{N}\) so that
\({\tau}_{\zeta}({\psi}^{-1} F) \leq D_\zeta\),
for every \(F \in S_\zeta[M]\), \(\max F \leq m_0^\zeta\).
Part 5. of Lemma \ref{L:21} now yields that 
\({\tau}_{\zeta}({\psi}^{-1} F) \leq
  (D_\zeta + k_\zeta)(\zeta + 1) + 1\),
for every \(F \in S_{\zeta}[M]\).
To complete the proof we need only take
\(D= \max \{(D_\zeta + k_\zeta)(\zeta + 1) + 1 : \zeta \leq \xi\}\). 
\end{proof}
\begin{Prop} \label{P:3}
Suppose \(L=(l_i)\), \(M=(m_i)\) belong to 
\([\mathbb{N}]\) and that \(\xi < \omega\).
Assume further that there exist a map
\(\psi \colon L \to M\)  
and an integer \(D > 0\)
so that \({\tau}_{\zeta}({\psi}^{-1} F) \leq D\),
for every \(F \in S_\zeta[M]\) and all
\(0 \leq \zeta \leq \xi\). Then there exist
integer constants \(E_\zeta\), \(0 \leq \zeta \leq \xi\),
so that
\({\tau}_{\zeta}({\phi}^{-1} F) \leq E_\zeta\),
for every \(F \in S_\zeta[M]\) and all
\(0 \leq \zeta \leq \xi\).
Here, \(\phi \colon L \to M\), is the natural bijection 
\(\phi(l_i) = m_i\).
\end{Prop}
\begin{proof}
If \(\zeta = 0\) the assertion is trivial 
(\(E_0=1\)). Suppose the assertion holds for some
\(\zeta \leq \xi - 1\). We will show that
\(E_{\zeta + 1}= [(\zeta + 1) E_\zeta + 1]
[(2D + 1)(\zeta + 2) + 1] \)
works for \(\zeta + 1\).
Let
\(F \in S_{\zeta + 1}[M]\), 
\(F= \{m_{i_1}, \dots , m_{i_p}\}\). Our 
hypothesis yields that \(\{l_{i_1}, \dots , l_{i_p}\}\)
is contained in the union of \(E_\zeta m_{i_1}\)
\(S_\zeta[L]\)-sets and so 
\({\tau}_{\zeta}(\{l_{i_1}, \dots , l_{i_p}\})
\leq (\zeta + 1) E_\zeta m_{i_1} +1\)
by part 5. of Lemma \ref{L:21}.
Choose \(q_1 \in \mathbb{N}\) so that the set
\(\{l_{i_1 + j}: 0 \leq j \leq q_1\}\)
is the union of exactly
\([(\zeta + 1) E_\zeta +1]
m_{i_1}\) successive
maximal \(S_\zeta[L]\)-sets. 

{\em CLAIM}. \({\tau}_{\zeta + 1}(
\{l_{i_1 + j}: 0 \leq j \leq q_1\}) \leq E_{\zeta + 1}\).

Once our claim is proven, we apply Lemmas \ref{L:2},
and \ref{L:21} (parts 4. and 1.) to conclude that
\({\tau}_{\zeta + 1}(\{l_{i_1 + j}: 0 \leq j \leq p-1 \})
\leq E_{\zeta +1}\) and
\({\tau}_{\zeta + 1}(\{l_{i_1}, \dots , l_{i_p}\})\)
\(\leq E_{\zeta +1}\).

To prove the claim we choose \(q < q_1\) so that
the set \(\{l_{i_1 + j}: 0 \leq j \leq q\}\) is
the union of exactly \(m_{i_1}\) successive, maximal
\(S_\zeta[L]\)-sets. Our task now is to show that
\({\tau}_{\zeta + 1}(
\{l_{i_1 + j}: 0 \leq j \leq q\}) \leq 
  (2D + 1)(\zeta + 2) +1\).
The claim will then follow by applying parts 4. and 2.
of Lemma \ref{L:21}.

We first observe that if \(0 \leq j_0 \leq q\) is chosen
so that \(l_{i_1 + j} \leq m_{i_1 +j}\), for all
\(j \leq j_0\), then 
\(\{l_{i_1 + j}: 0 \leq j \leq j_0\}\) 
is contained in the union of \(2D\)
\(S_{\zeta + 1}[L]\)-sets. 

Indeed, 
\(\{m_{i_1 + j}: 0 \leq j \leq j_0\}\) 
belongs to \(S_{\zeta + 1}[M]\), by part 3. of Lemma
\ref{L:21} and the fact that  
\({\tau}_{\zeta}(\{l_{i_1 + j}: 0 \leq j \leq j_0\}
\leq m_{i_1}\).
It follows now
that \(\Psi=\{ \psi(l_{i_1 + j}): \, 0 \leq j \leq j_0 , \,
\psi(l_{i_1 + j}) \geq m_{i_1} \}\) 
belongs to
\(S_{\zeta + 1}[M]\). 
To see this let \(\{m_{t_0} <  \cdots <  m_{t_k}\}\), where
\(k \leq j_0\) and \(i_1 \leq t_0\), be an enumeration of
\(\Psi\). Then \(m_{t_j} \geq m_{i_1 + j}\), for every
\(0 \leq j \leq k\). Since 
\(\{m_{i_1 + j}: 0 \leq j \leq k\}\) 
belongs to \(S_{\zeta + 1}[M]\) which is spreading,
we conclude that \(\Psi\) belongs to \(S_{\zeta + 1}[M]\).
Our hypothesis (for \(\zeta +1 \)) yields that
\({\psi}^{-1} (\Psi)= \{ l_{i_1 + j}: \, 0 \leq j \leq j_0 , \,
\psi(l_{i_1 + j}) \geq m_{i_1}\}\)
is contained in the union of
\(D\) \(S_{\zeta + 1}[L]\)-sets.
On the other hand, the cardinality of the set
\(\{ \psi(l_{i_1 + j}): \, 0 \leq j \leq j_0 , \,
\psi(l_{i_1 + j}) < m_{i_1} \}\) 
is at most \(i_1 -1\). 
Our hypothesis (for \(\zeta=0\)) now yields that 
the cardinality of the set
\(\{ l_{i_1 + j}: \, \psi(l_{i_1 + j}) < m_{i_1}, \,
0 \leq j \leq j_0\}\)
is at most \(D(i_1 -1)\). We deduce, since
\(l_{i_1} > i_1 -1\), that 
\(\{ l_{i_1 + j}: \, \psi(l_{i_1 + j}) < m_{i_1}, \,
0 \leq j \leq j_0\}\) is contained in the union
of \(D\) \(S_1[L]\)-sets. Hence,
\(\{l_{i_1 + j}: 0 \leq j \leq j_0\}\) 
is contained in the union of \(2D\)
\(S_{\zeta + 1}[L]\)-sets.

Next set \(j_1 = \min \{ j: \, 0 \leq j \leq q
\text{ and } l_{i_1 + j} > m_{i_1 + j} \}\).
If \(j_1\) does not exist, then
\(l_{i_1 + j} \leq m_{i_1 + j}\), for all
\(0 \leq j \leq q\). We obtain, by our previous observation
for \(j_0=q\), that 
\(\{l_{i_1 + j}: 0 \leq j \leq q\}\) 
is contained in the union of \(2D\)
\(S_{\zeta + 1}[L]\)-sets.

If \(j_1\) does exist, then
\(\{l_{i_1 + j}: 0 \leq j < j_1\}\) 
is contained in the union of \(2D\)
\(S_{\zeta + 1}[L]\)-sets.
Indeed, this is obvious if \(j_1=0\). If \(j_1 \geq 1\)
the assertion follows from our previous observation
by taking \(j_0=j_1 -1\).
Finally, \(\{l_{i_1 + j}: j_1 \leq j \leq q\}\) 
belongs to \(S_{\zeta + 1}[L]\), since
\(l_{i_1 + j_1} > m_{i_1 + j_1} \geq m_{i_1}\)
and  \(\{l_{i_1 + j}: j_1 \leq j \leq q\}\)
is contained in the union of \(m_{i_1}\)
successive maximal \(S_\zeta[L]\)-sets.
Thus, \(\{l_{i_1 + j}: 0 \leq j \leq q\}\) 
is contained in the union of \(2D + 1\)
\(S_{\zeta + 1}[L]\)-sets.

Concluding, in any case, the set
\(\{l_{i_1 + j}: 0 \leq j \leq q\}\) 
is contained in the union of \(2D + 1\)
\(S_{\zeta + 1}[L]\)-sets. 
Hence, applying
part 5. of Lemma \ref{L:21}, we obtain that
\({\tau}_{\zeta + 1}(
\{l_{i_1 + j}: 0 \leq j \leq q\}) \leq 
  (2D + 1)(\zeta + 2) +1\),
as desired. The proof of the proposition
is complete.
\end{proof}
\begin{Prop} \label{P:4}
Let \(\xi < \omega_1\) and \(L=(l_i)\), \(M=(m_i)\) be in
\([\mathbb{N}]\). Suppose that there exist \(\delta > 0\)
and a bounded linear operator
\(T \colon X_L^\xi \to X_M^\xi\)
such that \(\|T(e_l^\xi)\|_0 > \delta\),
for every \(l \in L\). Then there exist
a map \(\psi \colon L \to M\) and a bounded
linear operator \(R \colon X_L^\xi \to X_M^\xi\)
such that \(R(e_l^\xi) = e_{\psi(l)}^\xi\),
for every \(l \in L\).
\end{Prop}
\begin{proof}
Following \cite{LT}, given two infinite matrices \((a_{ij})\) and
\((d_{ij})\), we shall call \((d_{ij})\) a
block diagonal of \((a_{ij})\),
if there exist 
\((r_k)\), \((s_k)\), increasing sequences of positive
integers so that
\(d_{ij}=
\begin{cases}
a_{ij}, &\text{ if } (i,j) \in 
               \cup_{k=1}^{\infty}
           [r_k, r_{k+1}) \times [s_k, s_{k+1}) \\
 0, &\text{ otherwise. }
\end{cases}
\)
We can represent \(T\) as an infinite matrix
\((a_{ij})\). Then \(T(e_{l_i}^\xi) = 
\sum_{j=1}^{\infty} a_{ij} e_{m_j}^\xi\),
for every \(i \in \mathbb{N}\). Because
\(\|T(e_{l_i}^\xi)\|_0 > \delta \), for every 
\(i \in \mathbb{N}\) there exists \(j \in \mathbb{N}\)
such that \(|a_{ij}| > \delta\). We can thus define a map
\(\psi \colon L \to M\) so that if 
\(\psi(l_i)=m_j\), then \(|a_{ij}| > \delta\).
Observe that \({\psi}^{-1} \{m_j\}\) is finite, for all
\(j \in \mathbb{N}\), since \((T(e_{l_i}^\xi))\) is
weakly null in \(X_M^\xi\).
In particular, \(\psi(L) \in [M]\). Let
\((m_{k_j})_{j=1}^\infty\) be the increasing enumeration of
\(\psi(L)\). Given \(x= \sum_{i=1}^\infty \lambda_i e_{l_i}
\in X_L^\xi\), we set \(S(x)= \sum_{j=1}^\infty 
(\sum_{i=1}^\infty \lambda_i a_{i k_j}) e_{m_{k_j}}^\xi\).
It follows,
since \(T\) is bounded and \((e_n^\xi)\) is unconditional,
that \(S\) is a well defined bounded linear operator
from \(X_L^\xi\) into \(X_{\psi(L)}^\xi\).
Moreover, the matrix representation \((c_{ij})\)
of \(S\) with respect to the bases
\((e_{l_i}^\xi)\) and \((e_{m_{k_j}}^\xi)\)
is given by \(c_{ij}=a_{i k_j}\), for all positive integers
\(i, j\).

We next consider the matrix \((b_{ij})\) given by
\[ b_{ij}=
\begin{cases}
c_{ij}, &\text{ if } \psi(l_i)=m_{k_j} \\
0, &\text{ otherwise. }
\end{cases}
\]
Note that there exists a unique non-zero entry in every
row of the matrix \((b_{ij})\), while each column
contains only finitely many non-zero entries.
We can thus find \(p\), a permutation of
\(\mathbb{N}\), so that the matrix
\((b_{p(i) j})\) is a block diagonal of \((c_{p(i) j})\).
Since \((c_{p(i) j})\) represents the bounded linear
operator \(S \colon X_L^\xi \to X_{\psi(L)}^\xi\) 
with respect to the bases \((e_{l_{p(i)}}^\xi)\)
and \((e_{m_{k_j}}^\xi)\),
and \((b_{p(i) j})\) is a block
diagonal of \((c_{p(i) j})\), Proposition 1.c.8
of \cite{LT} yields that \((b_{p(i) j})\)
also represents a bounded linear operator from
\(X_L^\xi\) into \(X_{\psi(L)}^\xi\) 
with respect to the bases \((e_{l_{p(i)}}^\xi)\)
and \((e_{m_{k_j}}^\xi)\).
Consequently, \((b_{ij})\) represents a bounded
linear operator \(W \colon X_L^\xi \to X_{\psi(L)}^\xi\)
with respect to the bases \((e_{l_i}^\xi)\)
and \((e_{m_{k_j}}^\xi)\)
which evidently satisfies
\(W(e_{l_i}^\xi)= a_{i k_j} e_{\psi(l_i)}^\xi\) 
(where \(\psi(l_i)=m_{k_j}\)),
for all \(i \in \mathbb{N}\). 
Because \(|a_{i k_j}| > \delta\), if \(\psi(l_i)=m_{k_j}\),
and \((e_n^\xi)\) is unconditional, we obtain that
there exists a bounded
linear operator \(R \colon X_L^\xi \to X_M^\xi\)
such that \(R(e_l^\xi) = e_{\psi(l)}^\xi\),
for every \(l \in L\).
\end{proof}
We are now ready for the
\begin{proof}[Proof of Theorem \ref{T:1}]
\(3. \Rightarrow 1.\) and \(2. \Rightarrow 3.\)
are immediate. To prove that 1. implies 2.
we first apply Proposition \ref{P:4}
to obtain a map \(\psi \colon L \to M\) and 
a bounded
linear operator \(R \colon X_L^\xi \to X_M^\xi\)
such that \(R(e_l^\xi) = e_{\psi(l)}^\xi\),
for every \(l \in L\). Propositions \ref{P:2}
and \ref{P:3} will then yield a constant \(E > 0\)
such that \(\tau_\zeta(\phi^{-1} F) \leq E\), for
every \(F \in S_\zeta[M]\) and \(0 \leq \zeta \leq \xi\). 
(Where \(\phi \colon L \to M\)
is the natural bijection.)
The result now follows
from Lemma \ref{L:22}. 
\end{proof}
To obtain Corollary \ref{C:1} we shall need
the following
\begin{Lem} \label{L:28}
Let \(\xi < \omega\) and \(s=(u_n)\) be a bounded block basis of
\((e_n^\xi)\) such that \(\lim_n \|u_n\|_0 =0\). 
Then for every \(N \in [\mathbb{N}]\) and 
\(0 \leq \zeta \leq \xi\) there exists \(M \in [N]\)
so that \(\lim_n \|\zeta_n^M \cdot s\|_\zeta =0\).
(Given \(\mu= \sum_{i=1}^{\infty} a_i e_i \in c_{00}\),
we denote by \(\mu \cdot s\) the vector
\(\sum_{i=1}^{\infty} a_i u_i\) which of course belongs
to \(c_{00}\).)
\end{Lem}
\begin{proof}
If \(\zeta=0\), the assertion follows 
from the fact that \(\lim_n \|u_n\|_0 =0\).
Assume now that \(\zeta \leq \xi -1\) and that the
assertion holds for \(\zeta\). Let \(N \in [\mathbb{N}]\)
and \(\epsilon > 0\). We will find \(Q \in [N]\)
so that \(\|(\zeta +1)_1^Q \cdot s\|_{\zeta +1}
< \epsilon\). Once this is accomplished, we can choose
\((Q_i) \subset [N]\) so that  
\(\|(\zeta +1)_1^{Q_i} \cdot s\|_{\zeta +1}
< \epsilon_i\), where \(\lim_i \epsilon_i = 0\)
and \(F_1^{\zeta + 1}(Q_1) < F_1^{\zeta + 1}(Q_2) < \cdots \).
Letting \(M= \cup_{i=1}^{\infty} F_1^{\zeta + 1}(Q_i)\),
we obtain that \( (\zeta +1)_1^{Q_i} = (\zeta +1)_i^M\),
for all \(i \in \mathbb{N}\) and thus
\(\lim_i \|(\zeta +1)_i^M \cdot s\|_{\zeta +1} =0\).

We now pass to the construction of \(Q\).
By the induction hypothesis we can choose a sequence
\((P_i) \subset [N]\) satisfying the following properties:
\begin{enumerate}
\item \(F_1^\zeta(P_1) < F_1^\zeta(P_2) < \cdots \).
\item \(\min F_1^\zeta(P_1) > \frac{2 + 2b}{\epsilon}\),
where \(b\) is chosen so that \(\|u_n\|_\xi \leq b\),
for every \(n \in \mathbb{N}\).
\item \(\|\zeta_1^{P_i} \cdot s\|_\zeta < \frac{1}{2^i k_{i-1}}\),
for all \(i \geq 2\). Here we have set 
\(k_i = \max supp(\zeta_1^{P_i} \cdot s)\).
\end{enumerate}
Put \(Q= \cup_{i=1}^{\infty} F_1^\zeta(P_i)\).
We are going to show that 
\(\|\sum_{i=1}^n \zeta_i^Q \cdot s\|_{\zeta + 1}
\leq 2 + 2b\), for every \(n \in \mathbb{N}\).
Note that \(\zeta_i^Q = \zeta_1^{P_i}\) and
\(F_i^\zeta(Q) = F_1^\zeta(P_i)\), for all
\(i \in \mathbb{N}\).

Let \(G \in S_{\zeta + 1}\). Let also \(\{i_1, \dots , i_p\}\)
be an enumeration of \(\{i \leq n: supp(\zeta_i^Q \cdot s)
\cap G \ne \emptyset\}\). Choose \(l \leq \min G\)
and \(G_1 < \cdots < G_l\) in \(S_\zeta\) so that
\(G = \cup_{i=1}^l G_i\).
Then 
\(|(\zeta_{i_1}^Q \cdot s)(G)| \leq b\).
Further,
\begin{align}
\biggl| \sum_{t=2}^p \zeta_{i_t}^Q \cdot s (G) \biggr |
&= \biggl| \sum_{j=1}^l \sum_{t=2}^p  
    (\zeta_{i_t}^Q \cdot s) (G_j) \biggr| \leq
   \sum_{j=1}^l \sum_{t=2}^p |(\zeta_{i_t}^Q \cdot s) (G_j)|
\notag \\
&\leq \sum_{j=1}^l \sum_{t=2}^p \|\zeta_{i_t}^Q \cdot s\|_\zeta
\leq \sum_{j=1}^l \sum_{t=2}^p \frac{1}{2^{i_t}
 k_{i_t -1}} \notag \\
&\leq \frac{l}{k_{i_1}} \leq 1, \text{ since } 
l \leq \min G \leq k_{i_1}. \notag
\end{align}
Therefore, 
\(|\sum_{i=1}^n (\zeta_i^Q \cdot s) (G)|=
   |\sum_{t=1}^p (\zeta_{i_t}^Q \cdot s) (G)| \leq b + 1\).
It follows that \(\|\sum_{i=1}^n \zeta_i^Q \cdot s \|_{\zeta +1}
\leq 2 + 2b\). If we take \(n=\min F_1^\zeta (P_1)\),
we obtain that \(\|(\zeta + 1)_1^Q \cdot s \|_{\zeta + 1}
< \epsilon\). The proof of the lemma is now complete.
\end{proof}
\begin{proof}[Proof of Corollary \ref{C:1}]
Let \(T \colon X_L^\xi \to X_M^\xi\)
be an isomorphic embedding. We apply Theorem \ref{T:1} to
show that \((e_{l_n}^\xi)\) dominates \((e_{m_n}^\xi)\).
Indeed, we need only check that
\(\inf_{l \in L} \|T(e_l^\xi)\|_0 > 0\).
If that were not the case, let \((x_i)\) be
a subsequence of \((T(e_{l_i}^\xi))\) such that
\(\lim_i \|x_i\|_0 = 0\). By a standard perturbation
result we can assume, without loss of generality, that
for some block basis \((u_i)\) of \((e_{l_i}^\xi)\)
and a null sequence of positive scalars \((\epsilon_i)\)
we have that \(\|x_i - u_i\|_\xi < \epsilon_i\), for
all \(i \in \mathbb{N}\). It follows that also
\(\lim_i \|u_i\|_0 = 0\), and thus Lemma \ref{L:28}
yields \(N \in [\mathbb{N}]\) so that
\(\lim_i \|\xi_i^N \cdot s\|_\xi =0\), where \(s=(u_i)\).
But then \(\lim_i \|\xi_i^N \cdot x\|_\xi =0\) as well.
(\(x=(x_i)\)) This is a contradiction because
\((x_i)\) is equivalent to a subsequence of
\((e_{l_i}^\xi)\) and thus it is an \(\ell_1^\xi\)-spreading
model. 
Hence, \((e_{l_n}^\xi)\) dominates \((e_{m_n}^\xi)\) completing
the proof of part 1. 
Parts 2 and 3 are immediate consequences of Theorem
\ref{T:1}.
\end{proof}
We recall that a Banach space \(X\) is said to be 
{\em primary \/} if, for every bounded linear projection
\(P\) on \(X\), either \(PX\) or \((I-P)X\) is isomorphic
to \(X\).
\begin{Cor}
\(X_N^\xi\) is not primary, for every \(N \in [\mathbb{N}]\)
and all \(1 \leq \xi < \omega\).
\end{Cor}
\begin{proof}
We first let \(\mathcal{F}= \{(L,M) \in [N] \times [N]: 
L \cup M = N , \, L \cap M = \emptyset\}\).
\(\mathcal{F}\) is easily seen to be closed in
\([N] \times [N]\) and thus it is a Polish space.
We next set \(\mathcal{G}= \{ (L,M) \in \mathcal{F}:
   d_\xi(N,L)=d_\xi(N,M) = \infty \}\).
Arguing as we did in the proof of Lemma \ref{L:23}
we obtain that \(\mathcal{G}\) is a \(G_\delta\) dense
subset of \(\mathcal{F}\). If \((L,M) \in \mathcal{G}\)
then \(X_N^\xi = X_L^\xi \bigoplus X_M^\xi\).
However, Theorem \ref{T:1} implies that \(X_N^\xi\)
is not isomorphic to a subspace of either \(X_L^\xi\)
or \(X_M^\xi\).
\end{proof}
\section{Subspaces spanned by block bases}
In this section we investigate subspaces of \(X^\xi\)
spanned by block bases of \((e_n^\xi)\). We first show that
there exists a block basis of \((e_n^\xi)\) spanning a complemented
subspace of \(X^\xi\) which is not isomorphic to
\(X_M^\zeta\), for every \(M \in [\mathbb{N}]\) and
all \(0 \leq \zeta \leq \xi\).
\begin{Lem} \label{L:29}
Let \(x_1 < \cdots < x_p\) be a finite block basis of
\((e_n)\), the unit vector basis of \(c_{00}\). Let
also \(G_1 < \cdots < G_q\) be finite subsets of 
\(\mathbb{N}\) and \((a_i)_{i=1}^p\) be scalars.
Assume that there exists \(C > 0\) such that
\(|(\sum_{i \in I} a_i x_i)( \cup_{j \in J} G_j )| \leq C\),
whenever \(I \subset \{1, \dots , p\}\) and 
\(J \subset \{1, \dots , q\}\) satisfy one of the following
two conditions:
\begin{enumerate}
\item \(I = \cup_{j \in J} I_j\),
        \(I_{j_1} < I_{j_2}\) if \(j_1 < j_2\) and 
        \(I_j = \{ i \in I : suppx_i \cap G_j \ne \emptyset\}\),
 for all \(j \in J\).
\item \(J= \cup_{i \in I} J_i\),
        \(J_{i_1} < J_{i_2}\) if \(i_1 < i_2\) and 
        \(J_i = \{ j \in J : suppx_i \cap G_j \ne \emptyset\}\),
 for all \(i \in I\).
\end{enumerate}
Then \(|(\sum_{i=1}^p a_i x_i) (\cup_{j=1}^q G_j )| \leq 3 C\). 
\end{Lem}
\begin{proof}
Given \(j \leq q\), we let 
\(T_j = \{i \leq p: supp x_i \cap G_j \ne \emptyset\}\).
We also let 
\(J= \{ j \leq q: T_j \ne \emptyset \}\)
and \(J_1 = \{j \in J: |T_j | = 1 \}\).
Set \(J_2 = J \setminus J_1\). Given \(j \in J_2\)
we let \(s_j = \min T_j \) and \(t_j = \max T_j\).
We observe that \(s_{j_1} < t_{j_1} \leq s_{j_2}\), 
for every \(j_1 < j_2\) in \(J_2\).

Next, we define a map \(\sigma \colon J_1 \to \{1, \dots , p\}\)
so that \(\{\sigma (j) \} = T_j\), for every \(j \in J_1\). 
Note that \(\sigma (J_1)\) and \(J_1\) satisfy condition 2.
and therefore 
\[\biggl | \bigl (\sum_{i=1}^p a_i x_i \bigr )
  \bigl ( \cup_{j \in J_1} G_j  \bigr ) \biggr | =
  \biggl |\bigl (\sum_{i \in \sigma (J_1)} a_i x_i \bigr )
  \bigl ( \cup_{j \in J_1} G_j  \bigr ) \biggr | 
 \leq C. \]
Suppose now that \(J_2 = \{j_1, \dots , j_k\}\) and
put \(J_3 = \{j_r: r \leq k, r \text{ is odd } \}\)
and \(J_4 = \{j_r: r \leq k, r \text{ is even } \}\).
It follows that \(\cup_{j \in J_m} T_j\) and 
\(J_m\),
\(m \in \{3,4\}\), satisfy condition 1. and thus
\[\biggl | \bigl (\sum_{i=1}^p a_i x_i \bigr )
  \bigl ( \cup_{j \in J_m} G_j  \bigr ) \biggr | =
  \biggl | \bigl (\sum_{i \in \cup_{j \in J_m} T_j} a_i x_i \bigr )
  \bigl ( \cup_{j \in J_m} G_j  \bigr ) \biggr | 
 \leq C, \, m \in \{3,4\}.\]
Hence, \(| (\sum_{i=1}^p a_i x_i )
  ( \cup_{j \in J_2} G_j  ) | \leq 2C\).
The assertion follows since 
\[\bigl (\sum_{i=1}^p a_i x_i \bigr )
  \bigl ( \cup_{j=1}^q G_j \bigr )  =
 \bigl (\sum_{i=1}^p a_i x_i \bigr )
  \bigl ( \cup_{j \in J_1} G_j \bigr ) +
 \bigl (\sum_{i=1}^p a_i x_i \bigr )
  \bigl ( \cup_{j \in J_2} G_j \bigr ). \]
\end{proof} 
\begin{Lem} \label{L:30}
Let \(1 \leq \zeta \leq \xi < \omega\) and \((x_n)\) be a 
block basis of \((e_n^\xi)\) so that for some
\(b > 0\), \(\|x_n\|_\xi < b\), for every 
\(n \in \mathbb{N}\). Let also \(k_n = \max supp x_n\),
for every \(n \in \mathbb{N}\), and suppose that
\(\|x_n\|_{\zeta - 1} < \frac{1}{2^{k_{n-1}}}\), for every
\(n \geq 2\). Then 
\(| (\sum_{i=1}^n a_i x_i)(H) | \leq (2 + b) \max_{i \leq n}
|a_i|\), for every \(H \in S_\zeta\),
\(n \in \mathbb{N}\) and all scalar sequences
\((a_i)_{i=1}^n\).
\end{Lem}
\begin{proof}
Let \(H \in S_\zeta\) and put
\(i_0 = \min \{i \leq n: \, supp x_i \cap H \ne \emptyset \}\).
We may write \(H = \cup_{j=1}^r H_j \), where
\(r \leq \min H\) and \(H_1 < \cdots < H_r\) belong
to \(S_{\zeta - 1}\). Note that \(\min H \leq k_{i_0}\).
We also observe that
\(|x_i(H)| \leq r \|x_i\|_{\zeta - 1}\) and hence
\begin{align}
\bigl | \sum_{i= i_0 +1 }^n a_i x_i (H) \bigr | &\leq
\sum_{i= i_0 +1 }^n |a_i|r \|x_i\|_{\zeta - 1}
\notag \\
&\leq r(\max_{i \leq n} |a_i|) \sum_{i= i_0 +1 }^{\infty}
 \frac{1}{2^{k_{i-1}}} \leq
2 \max_{i \leq n} |a_i|. \notag
\end{align}
Finally, \(|x_{i_0}(H)| \leq \|x_{i_0}\|_\xi < b\) and thus

\(| (\sum_{i=1}^n a_i x_i)(H) | \leq (2 + b) \max_{i \leq n}
|a_i|\), as desired.
\end{proof}
Our next proposition is a partial generalization of
Lemma \ref{L:27}.
\begin{Prop} \label{P:5}
Let \(\xi < \omega\) and \((x_n)\) be a semi-normalized 
block basis of \((e_n^\xi)\).
Set \(\zeta = \min \{ \alpha \leq \xi: \inf_n \|x_n\|_\alpha
> 0\}\). Then there exists a subsequence of \((x_n)\)
which is equivalent to a subsequence of
\((e_n^{\xi - \zeta})\).
\end{Prop}
\begin{proof}
Choose \(\delta > 0\), \(b > 0\) so that
\(\delta < \|x_n\|_\zeta \) and
\(\|x_n\|_\xi < b\), for every \(n \in \mathbb{N}\).
Assume first that \(\zeta \geq 1\).
Then we choose inductively \(n_1 < n_2 < \cdots\)
so that 
\(\|x_{n_i}\|_{\zeta - 1} < \frac{1}{2^{k_{i-1}}}\), for every
\(i \geq 2\), where \(k_i = \max supp x_{n_i}\).
For every \(i \in \mathbb{N}\) we can find \(F_i \in S_\zeta\),
\(F_i \subset supp x_{n_i}\), so that 
\(|x_{n_i}|(F_i) > \delta\).
Put \(m_i = \min F_i\). We are going to show that
\((x_{n_i})\) is equivalent to 
\((e_{m_i}^{\xi - \zeta})\). To this end
let \(k \in \mathbb{N}\) and \((a_i)_{i=1}^k\)
be scalars. We first show that
\(\|\sum_{i=1}^k a_i e_{m_i}^{\xi - \zeta}\| \leq
{\delta}^{-1} \|\sum_{i=1}^k a_i x_{n_i} \|_\xi\).
Indeed, if \(G \subset \{m_1, \dots , m_k\}\)
belongs to \(S_{\xi - \zeta}\) then
set \(A= \{i \leq k: m_i \in G\}\). We have the
following estimate
\[ \sum_{i \in A} |a_i| \leq {\delta}^{-1}
    \sum_{i \in A} |a_i| |x_{n_i}|(F_i) \leq
   {\delta}^{-1} \|\sum_{i=1}^k a_i x_{n_i} \|_\xi\]
as \(\cup_{i \in A} F_i \in S_\xi\), by Lemma \ref{L:25}.

Next, let \(G \in S_\xi\). Lemma \ref{L:25} yields
\(G_1 < \cdots < G_q\) in
\(S_\zeta\) with
\(\{ \min G_j : j \leq q\}\)
belonging to \(S_{\xi - \zeta}\)
and so that \(G= \cup_{j=1}^q G_j\).
We shall apply Lemma \ref{L:29} in order to estimate
\(| (\sum_{i=1}^k a_i x_{n_i})(\cup_{j=1}^q G_j)|\).
Let \(I \subset \{1, \dots ,k\}\) and 
\(J \subset \{1, \dots , q\}\) satisfy condition 1.
of Lemma \ref{L:29}. 
Then \(I_j = \{i \in I : supp x_{n_i} \cap G_j \ne \emptyset\}\),
for every \(j \in J\). We choose \(i_j \in I_j\) such that
\(|a_{i_j}|= \max_{i \in I_j} |a_i|\), for every \(j \in J\).
Fix \(j_0 \in J\).
\begin{align}
\biggl | \bigl ( \sum_{i \in I_{j_0}} a_i x_{n_i} \bigr )
\bigl ( \cup_{j \in J} G_j \bigr ) \biggr | &=
\biggl | \sum_{i \in I_{j_0}} a_i x_{n_i} (G_{j_0}) \biggr |
\leq (2 + b) \max_{i \in I_{j_0}} |a_i| \notag \\
&= (2 + b) |a_{i_{j_0}}|,
\text{ by Lemma \ref{L:30} }. \notag
\end{align}
Hence 
\(| ( \sum_{i \in I} a_i x_{n_i} )
( \cup_{j \in J} G_j ) | \leq
(2 + b) \sum_{j \in J} |a_{i_j}|\).

Note also that \(\{m_{i_j} : j \in J \setminus \{\min J\}\}\)
belongs to \(S_{\xi - \zeta}\). This is so since
\(supp x_{n_i} \cap G_j \ne \emptyset\), whenever
\(i \in I_j\) and \(j \in J\), and thus
\(\min G_{j_1} < \min supp x_{n_i} \leq m_i\),
for every \(i \in I_{j_2}\) and \(j_1 < j_2\) in \(J\).
In particular, \(\min G_{j_1} < m_{i_{j_2}}\), 
when \(j_1 < j_2\) in \(J\).
Since \(S_{\xi - \zeta}\) is spreading we obtain that
\(\{m_{i_j} : j \in J \setminus \{\min J\}\}\)
belongs to \(S_{\xi - \zeta}\).
It follows now that \(\sum_{j \in J \setminus \{\min J \}} 
|a_{i_j}| \leq
\|\sum_{i=1}^k a_i e_{m_i}^{\xi - \zeta}\|\) and hence
\[\biggl | \bigl (\sum_{i \in I} a_i x_{n_i} \bigr )
\bigl (\cup_{j \in J} G_j \bigr ) \biggr |
\leq 2(2+b) 
\biggl \|\sum_{i=1}^k a_i e_{m_i}^{\xi - \zeta} \biggr \|. \]
We shall now assume that \(I \subset \{1, \dots , k\}\)
and \(J \subset \{1, \dots , q\}\)
satisfy condition 2. of Lemma \ref{L:29}.
Then \(J_i = \{j \in J: supp x_{n_i} \cap G_j \ne \emptyset\}\)
, for all \(i \in I\). An argument similar to that in the preceding 
paragraph, yields that
\(\{m_i : i \in I \setminus \{\min I\}\}\)
belongs to \(S_{\xi - \zeta}\).
It follows that \(\sum_{i \in I} |a_i| \leq
2 \|\sum_{i=1}^k a_i e_{m_i}^{\xi - \zeta}\|\).
Finally,
\begin{align}
\biggl | \bigl ( \sum_{i \in I} a_i x_{n_i} \bigr )
\bigl ( \cup_{j \in J} G_j \bigr ) \biggr | &=
\biggl | \sum_{i \in I} a_i x_{n_i} 
\bigl ( \cup_{j \in J_i} G_j \bigr ) \biggr |
\notag \\
&\leq b \sum_{i \in I} |a_i|,
\text{ as } \cup_{j \in J_i} G_j \in S_\xi,
\notag \\
&\leq 2b \biggl \|\sum_{i=1}^k a_i e_{m_i}^{\xi - \zeta} \biggr \|.
\notag
\end{align}
We deduce from Lemma \ref{L:29} that
\[\biggl | \bigl (\sum_{i=1}^k a_i x_{n_i} \bigr )
\bigl (\cup_{j=1}^q G_j \bigr ) \biggr | \leq
6(2+b) \biggl \|\sum_{i=1}^k a_i e_{m_i}^{\xi - \zeta} \biggr \|,\]
and hence \(\|\sum_{i=1}^k a_i x_{n_i}\|_\xi \leq 12(2+b) 
\|\sum_{i=1}^k a_i e_{m_i}^{\xi - \zeta}\|\). 

To complete the proof we need to consider the case
\(\zeta=0\). We now choose \(m_n \in supp x_n\)
such that \(|x_n|(\{m_n\}) > \delta\), for all
\(n \in \mathbb{N}\). We are going to show that
\((x_n)\) is equivalent to \((e_{m_n}^{\xi})\).
Arguing as we did in the case \(\zeta \geq 1\)
we obtain that  
\(\|\sum_{i=1}^k a_i e_{m_i}^{\xi}\| \leq
{\delta}^{-1} \|\sum_{i=1}^k a_i x_i \|_\xi\),
for every \(k \in \mathbb{N}\) and all scalar
sequences \((a_i)_{i=1}^k\).

Next let \(G \in S_\xi\) and put
\(I = \{i \leq k : suppx_i \cap G \ne \emptyset\}\).
Then 
\begin{align}
\biggl |\sum_{i=1}^k a_i x_i(G) \biggr | &\leq 
\sum_{i \in I} |a_i| | x_i(G)| \leq b \sum_{i \in I} |a_i|  \notag \\
&\leq
2b \biggl \|\sum_{i=1}^k a_i e_{m_i}^{\xi} \biggr \|, \text{ as }
\{m_i : i \in I \setminus \{\min I\}\} \in S_\xi.
\notag 
\end{align}
Hence \(\|\sum_{i=1}^k a_i x_i\|_\xi \leq
4b \|\sum_{i=1}^k a_i e_{m_i}^{\xi}\|\).
The proof of the proposition is now complete.
\end{proof}
As an immediate consequence of Proposition \ref{P:5}
we obtain
\begin{Cor} \label{C:4}
For every semi-normalized weakly null sequence in
\(X^\xi\), \(\xi < \omega\), there exist
\(\zeta \leq \xi\) and a subsequence which is
equivalent to a subsequence of \((e_n^\zeta)\).
\end{Cor}
\begin{Lem} \label{L:31}
Let \(1 \leq \xi < \omega\) and \((F_n)\) be a sequence
of successive members of \(S_\xi\) satisfying
the following
\begin{enumerate}
\item \(\bigl (\tau_{\xi - 1}(F_n) \bigr )\)
increases to \(\infty\).
\item \(\sup_n \frac{ \min F_n}{\tau_{\xi - 1}(F_{n+k})} > k\),
for every \(k \in \mathbb{N}\).
\end{enumerate}
Let \((u_n)\) be a convex block basis of
\((e_n^\xi)\) such that \(supp u_n = F_n\), for every
\(n \in \mathbb{N}\). Assume further that 
\(\sum_{n=1}^{\infty} \|u_n\|_{\xi -1} < \infty\).
Then the closed linear span of \((u_n)\) in \(X^\xi\)
is not isomorphic to \(X_M^\zeta\), for every
\(\zeta \leq \xi\) and \(M \in [\mathbb{N}]\).
\end{Lem}
\begin{proof}
\((u_n)\) is normalized in \(X^\xi\) since 
\(F_n \in S_\xi\), for every \(n \in \mathbb{N}\).
We let \(X\) denote the closed linear span of
\((u_n)\) in \(X^\xi\). Because 
\(\sum_{n=1}^{\infty} \|u_n\|_{\xi -1} < \infty\),
we deduce from Proposition \ref{P:5} that
every semi-normalized block basis of \((u_n)\)
admits a subsequence equivalent to the unit vector
basis of \(c_0\). Indeed, let \((v_n)\),
\(v_n = \sum_{i \in G_n} b_i u_i\), be a semi-normalized
block basis of \((u_n)\).
Note that \((b_n)\) is bounded since \((v_n)\) is.
But also, \(\lim_n \sum_{i \in G_n} \|u_i\|_{\xi -1} =0\),
since \(\sum_{n=1}^{\infty} \|u_n\|_{\xi -1} < \infty\).
Hence \(\lim_n \|v_n\|_{\xi -1} =0\) and therefore
Proposition \ref{P:5} (for \(\zeta=\xi\)) yields a subsequence
of \((v_n)\) equivalent to the unit vector basis of \(c_0\).

It follows that every semi-normalized weakly null sequence
in \(X\) admits a subsequence equivalent to the unit vector basis
of \(c_0\). That is, \(X\) has property \((S)\) \cite{cb}.
However, \(X_M^\zeta\) fails property \((S)\) when
\(\zeta \geq 1\) and \(M \in [\mathbb{N}]\).
Thus, \(X_M^\zeta\) is not isomorphic to a subspace of
\(X\) for every \(1 \leq \zeta \leq \xi\) and 
\(M \in [\mathbb{N}]\).

To complete the proof we show that \(X\) is
not isomorphic to \(c_0\). This is
accomplished by showing that for every \(k \in \mathbb{N}\)
there exists \(n \in \mathbb{N}\) so that
\((u_{n+i})_{i=1}^k\) is isometrically equivalent
to the unit vector basis of \({\ell}_1^k\).
In particular \(X\) contains uniformly complemented \({\ell}_1^k\)'s.
It is a well known fact that \(c_0\) fails this property.

We let \(k \in \mathbb{N}\) and choose according to 2.
\(n \in \mathbb{N}\) so that
\(\frac{ \min F_n}{\tau_{\xi - 1}(F_{n+k})} > k\).
Condition 1. now yields that 
\(\sum_{i=1}^k \tau_{\xi - 1}(F_{n+i}) < \min F_n\)
and thus \(\cup_{i=1}^k F_{n+i} \in S_\xi\). Hence
\[ \biggl \| \sum_{i=1}^k a_i u_{n+i} \biggr \|_\xi
\geq \sum_{i=1}^k |a_i|u_{n+i}(F_{n+i}) =
\sum_{i=1}^k |a_i|,\]
for every scalar sequence \((a_i)_{i=1}^k\).
Therefore \((u_{n+i})_{i=1}^k\) is isometrically equivalent
to the unit vector basis of \({\ell}_1^k\).
\end{proof}
\begin{Prop} \label{P:6}
Let \(1 \leq \xi < \omega\). There exists a normalized
convex block basis \((u_n)\) of \((e_n^\xi)\) so that letting
\(F_n = supp u_n\), for every \(n \in \mathbb{N}\), the following
are satisfied
\begin{enumerate}
\item \(\tau_{\xi -1}(F_n)=n^2\) and 
\(\min F_n > k(n+k)^2\), for every \(n\) and \(k\) in
\(\mathbb{N}\) such that \(k < n\).
\item \(X= [u_n : n \in \mathbb{N}]\) is not isomorphic to
\(X_M^\zeta\) for every \(\zeta \leq \xi\) and 
\(M \in [\mathbb{N}]\).
\item \(X\) is complemented in \(X^\xi\).
\end{enumerate}
\end{Prop}
\begin{proof}
We inductively choose a sequence of integer intervals
\((F_n)\) such that for every \(n \in \mathbb{N}\)
\[ \min F_n > \max \{k(n+k)^2: k < n\} \cup 
   \{\min F_{n-1}\} \text{ and }
    \tau_{\xi -1}(F_n)=n^2.\]
Put \(M_n = F_n \cup \{m \in \mathbb{N}: m > \max F_n\}\),
for every \(n \in \mathbb{N}\). We then define
\[ u_n = \frac{1}{n^2} \sum_{i=1}^{n^2}
         (\xi - 1)_i^{M_n}, \, n \in \mathbb{N}. \]
Condition 1. is an immediate consequence of the inductive
construction. This condition implies that in fact
\(F_n \in S_\xi\), for every \(n \in \mathbb{N}\) and thus
\((u_n)\) is indeed a normalized convex block basis of
\((e_n^\xi)\). We also obtain from Lemma \ref{L:3}
that \(\|u_n\|_{\xi - 1} \leq \frac{\xi}{n^2}\)
and so \(\sum_{n=1}^\infty \|u_n\|_{\xi - 1}
< \infty\). Hence
condition 2. holds
in view of Lemma \ref{L:31}. It remains to establish
that \(X\) is complemented in \(X^\xi\). To this end
we define a map \(P \colon c_{00} \to c_{00}\) by
\[ P(x)= \sum_{i=1}^{\infty} x(F_i) u_i, \text{ for all }
x \in c_{00}.\]
Clearly \(P\) is well defined and linear. It is also
clear that \(P(u_i)=u_i\), for every \(i \in \mathbb{N}\).
Our objective is to show that \(P\) is bounded with respect
to the \(\| \cdot \|_\xi \)-norm on \(c_{00}\), for then
\(P\) will extend to a bounded linear projection on
\(X^\xi\) with range equal to \(X\).
To achieve our goal it suffices to show that
if \(G \in S_\xi\) is maximal, then
\(( \sum_{i=1}^p x(F_i) u_i )
  (G) \leq 18 \xi \),
for every \(p \in \mathbb{N}\) and \(x \in c_{00}\), 
\(\|x\|_\xi \leq 1\),
with \(x(\{i\}) \geq 0\), \(i \in \mathbb{N}\).

According to condition 1. of our hypothesis,
for every \(i \in \mathbb{N}\) there exist
\(F_{i1} < \cdots < F_{ii^2}\) successive \(S_{\xi -1 }\)
sets so that \(F_i= \cup_{k=1}^{i^2} F_{ik}\)
, and $\{\min F_{ik} : k \leq i^2\} \in S_1$.
Next let \(q= \min G\) and choose 
\(G_1 < \cdots < G_q\) maximal members of \(S_{\xi -1 }\)
so that \(G= \cup_{j=1}^q G_j\). Of course 
\(\{\min G_j: j \leq q\}\) is maximal in \(S_1\). 
We shall apply Lemma \ref{L:29}.
Let \(I \subset \{1, \dots p\}\) and 
\(J \subset \{1, \dots q\}\) satisfy condition 1.
of Lemma \ref{L:29}. Recall that
\(I_j = \{i \in I: supp u_i \cap G_j \ne \emptyset\}\),
\(j \in J\).
For each \(j \in J\) we choose
\(i_j \in I_j\) and \(k_j \leq i_j^2\) such that
\(x(F_{i_j k_j})= \max_{k \leq i^2, \, i \in I_j} x(F_{ik})\).
We have the following estimate
\begin{align}
\biggl ( \sum_{i \in I} x(F_i) u_i \biggr )
\biggl ( \cup_{j \in J} G_j \biggr ) &=
\sum_{j \in J} \sum_{i \in I_j} x(F_i) u_i(G_j) \notag \\
&= \sum_{j \in J} \sum_{i \in I_j} 
    x(F_i) \frac{1}{i^2}
\sum_{k=1}^{i^2} (\xi - 1)_k^{M_i}(G_j) \notag \\
&\leq \sum_{j \in J} \sum_{i \in I_j} x(F_{i_j k_j})
      \sum_{k=1}^{i^2} (\xi - 1)_k^{M_i}(G_j), \notag \\ 
&\text{ since } F_i= \cup_{k=1}^{i^2} F_{ik}, \notag \\
&\leq \sum_{j \in J} x(F_{i_j k_j}) \sum_{i \in I_j}
      \sum_{k=1}^{i^2} (\xi - 1)_k^{M_i}(G_j)
\notag \\
&\leq \sum_{j \in J} x(F_{i_j k_j}) \xi, \text{ by Lemma }
\ref{L:3}, \notag \\
&\leq \xi x( \cup_{j \in J} F_{i_j k_j})
\leq 2 \xi.
\notag
\end{align}
The last inequality holds because \(\|x\|_\xi \leq 1\) and
\(\cup_{j \in J \setminus \{\min J\}} F_{i_j k_j}
\in S_\xi\). Indeed, \(\min G_{j_1} < \min F_{i_{j_2} k_{j_2}}\)
when \(j_1 < j_2\) in \(J\) and therefore, as
\(\{\min G_j: j \leq q\} \in S_1\),
\(\cup_{j \in J \setminus \{\min J\}} F_{i_j k_j}\) belongs
to \(S_\xi\) by Lemma \ref{L:25}.

Next assume that \(I \subset \{1, \dots , p\}\)
and \(J \subset \{1, \dots , q\}\) satisfy condition 2.
of Lemma \ref{L:29}. Then \(J_i = \{j \in J: supp u_i 
\cap G_j \ne \emptyset\}\), \(i \in I\).
We set \(H_i = \{ j \in J_i: G_j \subset suppu_i\}\),
\(i \in I\). 
Since $suppu_i = F_i$ is an interval, 
\(|J_i| \leq |H_i| + 2\), for all \(i \in I\).
Moreover, since 
each $G_j$ is a maximal $S_{\xi-1}$ set
and \(\tau_{\xi -1}(F_i)=i^2\), we have that
\(|H_i| \leq i^2\), for all \(i \in I\).
To estimate \(\sum_{i \in I} \sum_{j \in J_i \setminus H_i}
x(F_i) u_i(G_j)\), choose \(j_i \in J_i \setminus H_i\),
for every \(i \in I\) (we have assumed without loss of generality
that \(J_i \setminus H_i \ne \emptyset\)).
Then, the sets \(I\) and \(\{j_i : i \in I\}\) satisfy
condition 1. of Lemma \ref{L:29}. We deduce from our preceding
work that 
\[\sum_{i \in I} \sum_{j \in J_i \setminus H_i} 
  x(F_i) u_i(G_j) \leq 4 \xi,\]
as \(|J_i \setminus H_i| \leq 2\), for every \(i \in I\).
We next choose, for every \(i \in I\), 
\(R_i \subset \{1, \dots , i^2\}\) with \(|R_i|=|H_i|\)
and such that
\[\frac{1}{i^2} \sum_{k=1}^{i^2} x(F_{ik}) \leq
  \frac{1}{|H_i|} \sum_{k \in R_i} x(F_{ik}).\]
This choice is possible since \(|H_i| \leq i^2\).
(We make use of the following fact: Let
\((a_i)_{i=1}^n\) be scalars
with \(a_i \leq a_j\), \(i \leq j\), and let \(k < n\).
Then \(\frac{1}{n} \sum_{i=1}^n a_i \leq 
\frac{1}{n-k} \sum_{i=k+1}^n a_i\).) 
We now have that
\begin{align}
\sum_{i \in I} \sum_{j \in H_i} x(F_i) u_i(G_j)
&= \sum_{i \in I} \sum_{j \in H_i}
   \biggl [\sum_{k=1}^{i^2} x(F_{ik}) \biggr ]
   \biggl [ \frac{1}{i^2} \sum_{k=1}^{i^2} 
             (\xi - 1)_k^{M_i}(G_j) \biggr ]
\notag \\
&\leq \sum_{i \in I} \sum_{j \in H_i}
 \biggl [\frac{1}{i^2} \sum_{k=1}^{i^2} x(F_{ik}) \biggr ]
 \xi, \text{ by Lemma } \ref{L:3},
\notag \\
&\leq \xi \sum_{i \in I} \sum_{j \in H_i}
 \frac{1}{|H_i|} \sum_{k \in R_i} x(F_{ik}) \notag \\
&\leq \xi \sum_{i \in I} \sum_{k \in R_i} x(F_{ik})
\notag \\
&\leq \xi x(\cup_{i \in I, \, k \in R_i} F_{ik})
\leq 2 \xi. \notag
\end{align}
The latter inequality follows since \(\|x\|_\xi \leq 1\)
and \(\cup_{i \in I \setminus \{\min I\}, \, k \in R_i} F_{ik} 
\in S_\xi\). Indeed, the cardinality of the set
\(\{\min F_{ik}: k \in R_i, \, i \in I\}\) does not exceed
that of \(J\) since \(|R_i|=|H_i|\), for all 
\(i \in I\). It follows now, since
\(|J| \leq \min G_1\), that 
\(\{\min F_{ik}: k \in R_i, \, i \in I \setminus \{\min I\}\}\)
belongs to \(S_1\) and thus 
\(\cup_{i \in I, \, k \in R_i} F_{ik}\) is the union of
two members of \(S_\xi\). Concluding,
\begin{align}
\biggl ( \sum_{i \in I} x(F_i) u_i \biggr )
  ( \cup_{j \in J} G_j ) &=
\sum_{i \in I} \sum_{j \in J_i \setminus H_i}
x(F_i) u_i(G_j) + 
\sum_{i \in I} \sum_{j \in H_i}
x(F_i) u_i(G_j) \notag \\
&\leq 4 \xi + 2 \xi = 6 \xi. \notag
\end{align}
Lemma \ref{L:29} now implies that
\(( \sum_{i=1}^p x(F_i) u_i )
  (G) \leq 18 \xi \),
for every \(p \in \mathbb{N}\) and \(x \in c_{00}\), 
\(\|x\|_\xi \leq 1\),
with \(x(\{i\}) \geq 0\), \(i \in \mathbb{N}\).
It follows that \(\|P\| \leq 18 \xi\). The proof of
the proposition is now complete.
\end{proof}
\begin{Prop} \label{P:7}
Let \(1 \leq \xi < \omega\) and \((u_n)\)
be a block basis of \((e_n^\xi)\) satisfying 
the following
\begin{enumerate}
\item \(u_n = v_n + w_n\) with 
\(supp v_n \cap supp w_n = \emptyset\), \(n \in \mathbb{N}\).
\item \((w_n)\) is equivalent to the unit vector basis
of \(c_0\).
\item \(\lim_n \|v_n\|_\xi = 0\) yet
\(\sup_n \|\sum_{i=1}^n v_i \|_\xi =\infty\).
\end{enumerate}
Then there exists no projection from \(X^\xi\)
onto the closed linear span of \((u_n)\).
\end{Prop}
\begin{proof}
Let \(X\) denote the closed linear span of \((u_n)\)
in \(X^\xi\) and
assume that \(P \colon X^\xi \to X\) 
is a bounded linear projection.
Note that since \((e_n^\xi)\) is unconditional
our assumptions yield that \((u_n)\)
is semi-normalized in \(X^\xi\). Lemma 2.a.11
of \cite{LT} now yields that \((w_n)\) dominates
\((v_n)\) contradicting 3. as
\(\sup_n \|\sum_{i=1}^n w_i \|_\xi < \infty\). 
\end{proof}
It is easy to construct a normalized convex block
basis of \((e_n^\xi)\), \(\xi \geq 1\), satisfying
conditions 1.-3. Indeed, let \(M \in [\mathbb{N}]\),
\(M=(m_n)\), such that \(\sum_n \frac{1}{m_n} < \infty\).
Let \(q_n = \min F_n^\xi(M)\), \(n \in \mathbb{N}\)
(recall that \(F_n^\xi(M)= supp \xi_n^M\)). Because
\(\xi \geq 1\), \((e_{q_n}^\xi)\) is not dominated by the
unit vector basis of \(c_0\). It follows that there exists
a sequence of positive scalars \((a_n)\) such that
\(\lim_n a_n =0\) and \(\sup_n \|\sum_{i=1}^n a_i e_{q_i}^\xi \|
=\infty\). Set \(v_n = a_n e_{q_n}^\xi\) and
\(w_n = \frac{1-a_n}{1- \xi_n^M (q_n)} 
\sum_{i \in  F_n^\xi(M) \setminus \{q_n\}}
\xi_n^M (i) e_i^\xi\),
\(n \in \mathbb{N}\).
Finally, we let \(u_n = v_n + w_n\), \(n \in \mathbb{N}\).
Evidently, \((u_n)\) is a normalized convex block basis
of \((e_n^\xi)\) satisfying 1. and 3. It remains to show
that 2. holds. We observe that since
\(\sum_n \frac{1}{m_n} < \infty\) and 
\((\xi_n^M)\) is equivalent to the unit vector basis of
\(c_0\), then
letting \(x_n= 
\sum_{i \in  F_n^\xi(M) \setminus \{q_n\}}
\xi_n^M (i) e_i^\xi\), we have that
\(\sup_n \|\sum_{i=1}^n x_i \|_\xi < \infty\).
It follows that \((w_n)\) is equivalent to the 
unit vector basis of \(c_0\) as \(\lim_n a_n =0\)
and \(\lim_n \xi_n^M (q_n) = 0\).

\end{document}